\documentclass{article}[11pt]
\usepackage{amsthm,amsmath,a4wide}
\usepackage{amssymb}
\usepackage{amsfonts}
\usepackage[all,knot,tips]{xy}
\xyoption{arc}
\usepackage{mathrsfs}
\usepackage{bbm}
\usepackage{hyperref}
\input xypic
\usepackage[all,tips]{xy}
\usepackage{color}
\usepackage[titletoc]{appendix}

\numberwithin{equation}{section}

\newtheorem{thm}{Theorem}[section]
\newtheorem{lem}[thm]{Lemma}
\newtheorem{prop}[thm]{Proposition}
\newtheorem{cor}[thm]{Corollary}
\newtheorem{conj}[thm]{Conjecture}

\theoremstyle{definition}
\newtheorem{defin}[thm]{Definition}

\theoremstyle{remark}
\newtheorem{remark}[thm]{Remark}
\newtheorem{example}[thm]{Example}

\newcommand{\bth}{\begin{thm}}
\renewcommand{\eth}{\end{thm}}
\newcommand{\bpr}{\begin{prop}}
\newcommand{\epr}{\end{prop}}
\newcommand{\ble}{\begin{lem}}
\newcommand{\ele}{\end{lem}}
\newcommand{\bco}{\begin{cor}}
\newcommand{\eco}{\end{cor}}
\newcommand{\bde}{\begin{defin}}
\newcommand{\ede}{\end{defin}}
\newcommand{\bex}{\begin{example}}
\newcommand{\eex}{\end{example}}
\newcommand{\bre}{\begin{remark}}
\newcommand{\ere}{\end{remark}}
\newcommand{\bcj}{\begin{conj}}
\newcommand{\ecj}{\end{conj}}

\newcommand{\be}{\begin{equation}}
\newcommand{\ee}{\end{equation}}
\newcommand{\beq}{\begin{equation}}
\newcommand{\eeq}{\end{equation}}
\newcommand{\ve}{{\varepsilon}}
\newcommand{\ot}{{\otimes}}
\newcommand{\op}{{\oplus}}
\newcommand{\lb}{\label}
\newcommand{\bpf}{\begin{proof}}
\newcommand{\epf}{\end{proof}}

\newcommand{\uu}{{\mathfrak u}}

\newcommand{\s}{{\mathfrak s}}

\newcommand{\E}{{\cal E}}

\newcommand{\T}{{\cal T}}
\newcommand{\C}{{\cal C}}
\newcommand{\Z}{{\cal Z}}
\newcommand{\D}{{\cal D}}

\newcommand{\bZ}{{\mathbb Z}}
\newcommand{\bC}{{\mathbb C}}

\newcommand{\V}{{\cal V}}

\newcommand{\Rep}{{\cal R}{\it ep}}

\newcommand{\Funct}{{\cal F}{\it unct}}

\newcommand{\Aut}{{\cal A}{\it ut}}
\newcommand{\kk}{\mathbbm{k}}

\newcommand{\TL}{\mathcal{T\hspace{-1.5pt}L}}

\newcommand{\MF}{\mathrm{MF}}
\newcommand{\ZMF}{\mathrm{ZMF}}
\newcommand{\HMF}{\mathrm{HMF}}

\newcommand{\MFbi}{\mathrm{MF}_\mathrm{bi}}
\newcommand{\ZMFbi}{\mathrm{ZMF}_\mathrm{bi}}
\newcommand{\HMFbi}{\mathrm{HMF}_\mathrm{bi}}

\newcommand{\MFgr}{\mathrm{MF}^\mathrm{gr}}
\newcommand{\ZMFgr}{\mathrm{ZMF}^\mathrm{gr}}
\newcommand{\HMFgr}{\mathrm{HMF}^\mathrm{gr}}

\newcommand{\ZMFbigr}{\mathrm{ZMF}^\mathrm{gr}_\mathrm{bi}}
\newcommand{\HMFbigr}{\mathrm{HMF}^\mathrm{gr}_\mathrm{bi}}

\newcommand{\Pd}{\mathcal{P}_d}
\newcommand{\Pdgr}{\mathcal{P}_d^\mathrm{gr}}

\newcommand{\Ad}{Ad}
\newcommand{\id}{1}
\newcommand{\Id}{I\hspace{-.5pt}d}

\begin{document}

\vspace*{3em}

\begin{center}
{\LARGE $N=2$ minimal conformal field theories \\[.3em] and matrix bifactorisations of $x^d$}

\vspace{3em}

{\large 
Alexei Davydov$^{a}$,\, 
Ana Ros Camacho$^{b}$,\, 
Ingo Runkel$^{b}$ 
~\footnote{Emails: 
{\tt alexei1davydov@gmail.com}, 
{\tt anaroscamacho@gmail.com},
{\tt ingo.runkel@uni-hamburg.de}
}}
\\[2em]
\it$^a$ 
Department of Mathematics, Ohio University,\\
 Athens, OH 45701, USA
\\[1em]
$^b$ Fachbereich Mathematik, Universit\"at Hamburg,\\
Bundesstra\ss e 55, 20146 Hamburg, Germany
\end{center}

\vspace{3em}

\begin{abstract}
We prove a tensor equivalence between full subcategories of a) graded matrix factorisations of the potential $x^d-y^d$ and b) representations of the $N=2$ minimal super vertex operator algebra at central charge $3-6/d$, where $d$ is odd. The subcategories are given by a) permutation-type matrix factorisations with consecutive index sets, and b) Neveu-Schwarz-type representations.
The physical motivation for this result is the Landau-Ginzburg / conformal field theory correspondence, where it amounts to the equivalence of a subset of defects on both sides of the correspondence. Our work builds on results by Brunner and Roggenkamp
\cite{brunnerrogg},
where an isomorphism of fusion rules was established.
\end{abstract}

\newpage

\tableofcontents

\section{Introduction}

In this paper we will establish a tensor equivalence between certain categories of matrix bifactorisations and of representations of $N=2$ minimal super vertex operator algebras. Physically, this amounts to comparing the behaviour of a subset of defects at two ends of a renormalisation group flow. In this introductory section we will briefly review the physical motivation and provide some context for our result. The main body of the paper is purely mathematical and makes no further reference to the physical motivation.

\medskip

The main result of this paper can be seen as an instance of the so-called Landau-Ginzburg/ conformal field theory correspondence, which amounts to the following physical considerations. One starts from a family of two-dimensional quantum field theories, called $N=2$ supersymmetric Landau-Ginzburg models with target space $\bC^n$ and superpotential $W \in \bC[x_1,\dots,x_n]$ (see e.g.\ \cite{Hori:2003ic}). These theories are not conformally invariant and hence each such theory actually provides a one-parameter family of theories via the renormalisation group flow. Following the flow towards the short-distance behaviour (the UV theory) one reaches the free $N=2$ supersymmetric theory with target $\bC^n$. Following the flow towards the long-distance behaviour (the IR theory) results in an a priori unknown and typically
non-free $N=2$ superconformal field theory. 

By Zamolodchikov's $c$-theorem \cite{Zamolodchikov:1986gt}, the easiest statement to make about the IR theory is that its Virasoro central charge is strictly less than $3n$. Using quantities which stay invariant along the flow, one can deduce various other properties of the IR theory in terms of the initial data $n$ and $W$. For the purpose of this paper, let us single out three of these, in historical order:
\begin{enumerate}
\item  {\em Algebras of chiral primary fields} \cite{m1989,vw,Lerche:1989uy}: In the space of bulk fields one considers the subspace of chiral primaries as determined by $N=2$ supersymmetry. These fields have regular operator product expansion, resulting in the structure of an associative unital algebra over $\bC$ on this subspace, called the (c,c)-chiral ring. In terms of our initial data, it is given by the Jacobi ring
$$
	Jac(W) = \bC[x_1,\dots,x_n] \,/\, \big\langle \tfrac{\partial}{\partial x_1} W , \cdots ,  \tfrac{\partial}{\partial x_n} W \big\rangle \ .
$$
\item  {\em Categories of boundary conditions preserving B-type supersymmetry}
\cite{kl0210,bhls0305,hl0404}:
This is a $\bC$-linear category whose objects are boundary conditions that preserve the B-type subalgebra of the supersymmetry algebra. The morphisms are given by the $\bC$-linear subspace of chiral primaries amongst all boundary (changing) fields. The composition of morphisms is again obtained from the operator product expansion. In terms of our initial data, the category of boundary conditions is the homotopy category of matrix factorisations of $W$ over the algebra $\bC[x_1,\dots,x_n]$,
$$
	\HMF_{\bC[x_1,\dots,x_n],W} \ ,
$$
whose definition we recall in Section \ref{sec:mf}.
\item  {\em Tensor categories of defect conditions preserving B-type supersymmetry} \cite{brunnerrogg,cr1}:  This is a $\bC$-linear tensor category whose objects are defect conditions compatible with B-type supersymmetry. Morphisms and composition are defined as for boundary conditions. The zero distance limit of two defect lines defines the so-called fusion of defects, providing a tensor product on the category of defect conditions. In terms of our initial data, the tensor category of defects is a ``bimodule version'' of the above category,
$$
	\HMF_{\text{bi};\bC[x_1,\dots,x_n],W} \ ,
$$
see Section \ref{sec:mf}.
\end{enumerate}

If one has some independent access to  the above quantities in a candidate IR theory, one can try to compare them to the Landau-Ginzburg results given above. 
For example, if the candidate IR theory is a rational conformal field theory, the representation theory of vertex operator algebras provides such an alternative approach, leading to surprising mathematical statements. 

In this paper we are concerned with the third of the above invariants, and in this case it is convenient to use a graded variant of the above construction. The grading is provided by the so-called R-charge. On the Landau-Ginzburg side, $W$ is quasi-homogeneous of total degree 2, and the degrees $|x_i|$ of its variables now form part of our initial data. We then restrict to the R-charge zero sector in the invariants 1--3 above. Invariant 1, the chiral ring, then becomes trivial -- its R-charge zero subalgebra is just $\bC$. But invariants 2 and 3 remain interesting, see Section \ref{sec:mf} for the definitions of $\HMFgr_{\bC[x_1,\dots,x_n],W}$  and $\HMFgr_{\text{bi};\bC[x_1,\dots,x_n],W}$.

On the conformal field theory side, there is an elegant description of boundary conditions and defects in the case of bosonic, non-supersymmetric theories \cite{tft1,Frohlich:2006ch}\footnote{
Unfortunately, the corresponding description for rational $N=2$ superconformal theories has to date not been worked out. But one may reasonably expect that the result will be similar. For the sake of this introduction, we use the bosonic description as a placeholder for the yet-to-be-given supersymmetric variant. We also note that for unitary $N=2$ theories, conformal weight zero implies R-charge zero.
\label{fn:no-susy}
}. The initial data in this case is a rational vertex operator algebra $V$, together with a $\bC$-linear category $\mathcal{B}$, which is in addition a module category over $\Rep(V)$.
\begin{enumerate}  \setcounter{enumi}{1}
\item
Boundary conditions preserving $V$ are described by $\mathcal{B}$ itself. For $M,N \in \mathcal{B}$, the morphism space $\mathcal{B}(M,N)$ is the space of conformal weight zero boundary fields changing $M$ to $N$.
\item
Defect conditions transparent to the holomorphic and anti-holomorphic copy of $V$ are described by the tensor category $\mathcal{E}\hspace{-.8pt}nd_{\Rep(V)}(\mathcal{B})$ of module-category endofunctors of $\mathcal{B}$. Natural transformations of module functors describe the conformal weight zero defect (changing) fields.
\end{enumerate}
An important example (and in fact the example relevant to this paper) is provided by choosing $\mathcal{B}=\Rep(V)$ (as module category over itself) which entails  $\mathcal{E}\hspace{-.8pt}nd_{\Rep(V)}(\mathcal{B}) \simeq \Rep(V)$ (as tensor categories).

\medskip

There are two reasons why one should not expect an equivalence between the two descriptions of invariants 2 and 3 given above. Firstly, the boundary (defect) conditions above have the extra requirement of compatibility with $V$ (respectively $V \otimes_{\bC} V$), and the renormalisation group flow end points of Landau-Ginzburg boundary (defect) conditions may or may not satisfy this requirement. Secondly, not all boundary (defect) conditions of the IR theory may arise as end points of renormalisation group flows. The prediction, therefore, is that (up to footnote \ref{fn:no-susy}):
\begin{enumerate}  \setcounter{enumi}{1}
\item
a full subcategory of $\HMFgr_{\bC[x_1,\dots,x_n],W}$ is equivalent, as a $\bC$-linear category, to a full subcategory of $\mathcal{B}$;
\item
a full tensor subcategory of $\HMFgr_{\text{bi};\bC[x_1,\dots,x_n],W}$ is equivalent, as a $\bC$-linear tensor category, to a full tensor subcategory of $\mathcal{E}\hspace{-.8pt}nd_{\Rep(V)}(\mathcal{B})$.
\end{enumerate}
The point we wish to make in treating invariants 2 and 3 alongside each other is that {\em invariant 3 is much stronger} as it compares $\bC$-linear {\em tensor} categories. Still, there are suprisingly few examples where  even only a correspondence of some objects in $\HMFbigr$ and $\mathcal{E}\hspace{-.8pt}nd_{\Rep(V)}(\mathcal{B})$ is provided \cite{brunnerrogg,Behr:2014bta}. And, prior to the present work, there was no example in which a tensor equivalence has been established (beyond group-like subcategories, cf.\ \cite{cr1}). 

\medskip

Let us now describe the mathematical contents of this paper in more detail. On the Landau-Ginzburg side, we consider the case that $W$ depends only on a single variable $x$ and is given by $W = x^d$. In $\HMFgr_{\text{bi};\bC[x],x^d}$ we select the full tensor subcategory $\Pdgr$ which consists of so-called permutation-type matrix factorisations which have consecutive index sets,
see Section \ref{sec:mf}. On the conformal field theory side, we take the bosonic part of the
$N=2$ minimal super vertex operator algebra $V(N{=}2,d)$ and consider the full tensor subcategory $\C(N{=}2,d)_{NS}$ of its NS-type representations. Our main result is

\medskip\noindent
{\bf Theorem \ref{thm:main}.} 
{\em For $d$ odd, there is a tensor equivalence $\Pdgr \simeq \C(N{=}2,d)_{NS}$.}

\medskip

Our work is based on \cite{brunnerrogg}, where (for all $d$) the existence of a multiplicative equivalence (that is, a functor for which $F(M \otimes N) \simeq F(M) \otimes F(N)$, but without a coherence condition on the isomorphisms) is established. The missing piece provided by the above theorem is the comparison of associators. For some specific triples of objects (but for all $d$), this comparison of associators was already carried out in \cite{cr1}.

The proof of Theorem \ref{thm:main} works by first establishing a universal property for $\C(N{=}2,d)_{NS}$, that is we describe tensor functors out of $\C(N{=}2,d)_{NS}$.
We do this by means of universal properties of Temperley-Lieb categories and products with pointed categories (Section \ref{sec:N=2-rep}). This description requires $d$ to be odd. We then use the universal property to obtain a tensor functor into $\Pdgr$ and use a semi-simplicity argument to show that it is an equivalence (Section \ref{sec:mf}). 

\medskip

Returning for a moment to the general discussion of invariants 2 and 3 above, we wish to point out that there is currently no general mechanism known to find which potentials $W$ correspond to which pairs $(V,\mathcal{B})$, nor
a criterion to single out the relevant subcategories. It would of course be highly desirable to prove the equivalences in 2 and 3 without working out both sides explicitly first, but this seems currently out of reach.

\subsection*{Acknowledgements}
We would like to thank
	Hanno Becker,
	Nils Carqueville,
	Scott Morrison,
	Jeffrey Morton,
	Daniel Murfet
and
	Daniel Roggenkamp
for helpful discussions, and Nils Carqueville for many useful comments on a draft version of this paper.
AD would like to thank the Max Planck Institut f\"ur Mathematik (Bonn) for hospitality and excellent working conditions during the final stage of the work on the paper. AD also thanks the Department of Mathematics of Hamburg University for hospitality during visits in 2013 and 2014. These visits were partially supported by the Research Training Group 1670 ``Mathematics inspired by string theory and quantum field theory'' of the Deutsche Forschungsgemeinschaft.
ARC is supported by the Research Training Group 1670.
IR is supported in part by the DFG funded Collaborative Research Center 676 ``Particles, Strings, and the Early Universe''.

\subsection*{Notation}

Let $\kk$ be an algebraically closed field (which can be assumed to be the field $\bC$ of complex numbers).
We call a category $\C$ {\em tensor} if it is an additive
$\kk$-linear monoidal category such that the tensor product is $\kk$-linear in both arguments. 
A monoidal functor between tensor categories is {\em tensor} if it is $\kk$-linear. 
By {\em fusion} category we mean a tensor category which is semi-simple, with finite dimensional morphism spaces, finitely many isomorphism classes of
simple objects, and with simple unit object.

\section[Categories of representations for $\boldsymbol{N{=}2}$ minimal super vertex operator algebras]{Categories of representations for $\boldsymbol{N{=}2}$ minimal\newline super vertex operator algebras}\label{sec:N=2-rep}

\subsection[Representations of $N=2$ minimal super vertex operator algebras]{Representations of $\boldsymbol{N=2}$ minimal super vertex operator algebras}

Let $V(N{=}2,d)$ be the super vertex operator algebra of the $N=2$ minimal model of central charge $c=\frac{3(d-2)}{d}$, where $d\in\bZ_{\ge 2}$, see \cite{ad1} and e.g.\ \cite{dv,Eholzer:1996zi,Adamovic:1998} for more on $N=2$ superconformal algebras.
Its bosonic part $V(N{=}2,d)_0$ can be identified with the coset $(\widehat{\s\uu(2)}_{d-2}\oplus \widehat{\uu(1)}_4)/\widehat{\uu(1)}_{2d}$ \cite{dv} (see \cite{Carpi:2012va} for a proof in the framework of conformal nets).

Accordingly, the category $\C(N{=}2,d)$ of representations of $V(N{=}2,d)_0$ can be realised as the category of local modules over a commutative algebra $A$ in the product
\begin{align}
\E &= \Rep(\widehat{\s\uu(2)}_{d-2})\  \boxtimes\ \overline{\Rep(\widehat{\uu(1)}_{2d})} \boxtimes\ \Rep(\widehat{\uu(1)}_4) \nonumber \\ 
& = \C(\s\uu(2),{d{-}2})  \boxtimes\ \C(\bZ_{2d},q^{-1}_{2d}) \boxtimes\ \C(\bZ_4,q_4)\ ,
\label{eq:E-def}
\end{align}
see \cite{Frohlich:2003hg}.
Here, for a ribbon category $\C$ the notation $\overline\C$ stands for the tensor category $\C$ with the opposite braiding and ribbon twist. 
The category $\C(\s\uu(2),d{-}2)= \Rep(\widehat{\s\uu(2)}_{d-2})$ is the category of integrable highest weight representations of the affine $\s\uu(2)$ at level $d-2$.
Its simple objects $[l]$ are labelled by $l=0,...,d-2$ and have lowest conformal weight $h_l=\frac{l(l+2)}{4d}$. Their dimensions are $\mathrm{dim}[l] = \frac{\eta^{l+1}-\eta^{-l-1}}{\eta-\eta^{-1}}$
with $\eta = e^{ 2 \pi i /d }$ and their ribbon twists are $\theta_l = e^{2\pi i h_l} \, \id_{[l]}$. The fusion rule of $\C(\s\uu(2),d{-}2)$ is 
$$[k] \otimes [l] ~\simeq~ \bigoplus_{m=|k-l|~\mathrm{step}\,2}^{\min(k+l,2d-4-k-l)} [m] \ .$$
The category $\Rep(\widehat{\uu(1)}_{2d})$ of representations of the vertex operator algebra for $\uu(1)$, rationally extended by two fields of weight $d$, is a pointed fusion category (a fusion category with a group fusion rule) with group $G$ of isomorphism classes of simple objects given by $\bZ_{2d}$. Braided monoidal structures on pointed fusion categories require $G$ to be abelian and
are classified by quadratic functions $q:G\to\bC^*$ \cite{js}. The ribbon twist of $\C(G,q)$ is $\theta_a = q(a)\,\id$. The $q_m$ appearing in \eqref{eq:E-def} are defined as
$q_m:\bZ_m\to\bC^*$ with $q_m(r) = e^\frac{\pi ir^2}{m}$ and $m$ even.

We can label simple objects of $\E$ by $[l,r,s]$, where $l \in \{0,...,d-2\}$, $r \in \bZ_{2d}$ and $s \in \bZ_4$. 
The ribbon twist for $\E$ is given by $\theta_{[l,r,s]} = e^{2 \pi i h_{l,r,s}} \,\id$ with
$$h_{l,r,s} \equiv \frac{l(l+2)}{4d} + \frac{s^2}{8} - \frac{r^2}{4d} \mod \bZ\ .$$ 
The underlying object of the algebra $A$ 
in the product \eqref{eq:E-def}
is $[0,0,0]\oplus[d{-}2,d,2]$. Note that $[d{-}2,d,2]$ is an invertible object of order 2 and ribbon twist 1, so that $[0,0,0]\oplus[d{-}2,d,2]$ has a uniquely defined commutative separable algebra structure. The tensor product with $[d{-}2,d,2]$ has the form
$$[d{-}2,d,2] \otimes [l,r,s] ~\simeq~ [d{-}2{-}l,r{+}d,s{+}2].$$
In particular no simple objects are fixed by tensoring with $[d{-}2,d,2]$ and hence all simple $A$-modules are free:
\beq\lb{moa}
A\otimes[l,r,s] ~\simeq~ A\otimes [d{-}2{-}l,r{+}d,s{+}2] ~\simeq~ [l,r,s]\ \op\ [d{-}2{-}l,r{+}d,s{+}2]\ .
\eeq

Recall that a simple $A$-module is {\em local} if all its simple constituents have the same 
ribbon twist (see \cite{Pareigis,Kirillov:2001ti} and \cite[Cor.\,3.18]{Frohlich:2003hm}).
Thus local $A$-modules correspond to $[l,r,s]$ with even $l+r+s$:
$$h_{d-2-l,r+d,s+2} - h_{l,r,s} = \frac{(d{-}2{-}l)(d{-}l)-l(l{+}2)}{4d} + \frac{(s{+}2)^2-s^2}{8} - \frac{(r{+}d)^2-r^2}{4d} 
=\frac{s-l-r}{2}\ .$$

The fermionic part $V(N{=}2,d)_1$ of $V(N{=}2,d)$ corresponds to the $A$-module 
$$A \otimes [0,0,2] ~\simeq~ [0,0,2]\ \op\ [d-2,d,0]$$ 
so that the simple objects of the NS (R) sector of $\C(N{=}2,d)$ are $A\ot[l,r,s]$ with even (odd) $s$:
$$h_{l,r,s+2} - h_{l,r,s} - h_{0,0,2} = \frac{(s+2)^2-s^2-4}{8} = \frac{s}{2}\ .$$
Denote by $\C(N{=}2,d)_{NS}$ the full subcategory of $\C(N{=}2,d)$ consisting of NS objects, i.e.\ with simple objects of the form $A\ot[l,r,s]$ with even $s$. 
By \eqref{moa} any simple object in $\C(N{=}2,d)_{NS}$ can be written as 
\be\label{eq:[l,r]-def} 
[l,r] := A\ot[l,r,0] 
\quad \text{with}
\quad l \in \{0,1,\dots,d-2\} ~,~~ r \in \bZ_{2d} ~,~~ l+r \text{ even} \ .
\ee

\subsection{The structure of $\C(N{=}2,d)_{NS}$ for odd $d$}

Note that direct sums of objects $[l,r,s]$ with even $l+r+s$ form a ribbon fusion subcategory $\E_{even}$ of $\E$.
It can be characterised as the M\"uger centraliser of $[d{-}2,d,2]$ in $\E$. Recall that the M\"uger centraliser of a subcategory $\D\subset\C$ in a ribbon fusion category is $\{X\in\C\ |\ \theta_{X\ot Y} = \theta_X\ot\theta_Y,\ \forall Y\in\D\}$ \cite{mu}. 

The induction functor $A\ot-:\E\to {_A\E}$ is a faithful tensor functor. Its restriction to $\E_{even}$ is in addition ribbon, so that
$$\E_{even}~ \xrightarrow{A \ot -} ~ {_A\E}_{even} = {_A\E}^{loc} = \C(N{=}2,d)$$
is a faithful ribbon tensor functor.
For odd $d$ 
the object $[1,d,0]$ lies in $\E_{even}$ and tensor generates a subcategory 
of $\E_{even}$ with simple objects $[l,dl,0],\ l=0,...,d-2$ and the fusion with $[1,d,0]$ given by
\be\label{eq:E-fusion}
[1,d,0]\otimes [l,dl,0] \simeq 
\begin{cases}  [l{-}1,d(l{-}1),0]\ \op\ [l{+}1,d(l{+}1),0]  &; ~~ 1\leq l < d-2 \\  { [d{-}3,d(d{-}3),0] }
 & ; ~~ l = d-2\end{cases}
 \ee
Since the last entry in $[l,dl,0]$ is zero,
the restriction of the induction functor $A\ot-$ to this subcategory is fully faithful. Denote by $\T$ its image in $\C(N{=}2,d)$.

The invertible object $[0,2,0]$ belongs to the M\"uger centraliser of $[1,d,0]$ in $\E_{even}$:
$$\exp 2 \pi i \big(h_{1,d+2,0} - h_{1,d,0} - h_{0,2,0}\big) = \exp 2 \pi i\big(\tfrac{(d+2)^2-d^2-4}{4d}\big) = 1 \ .$$
It tensor generates a pointed subcategory $\V$ in $\E_{even}$ equivalent to $\C(\bZ_d,q^{-2}_d)$. The restriction of the induction functor $A\ot-$ to this subcategory is fully faithful.

For $d$ odd, $[1,d] \in \C(N{=}2,d)_{NS}$ and it is straightforward to see that $\C(N{=}2,d)_{NS}$ is tensor generated by $[1,d]$ and $[0,2]$ (recall the notation \eqref{eq:[l,r]-def}). Furthermore, the intersection of the subcategories tensor generated by $[1,d]$ and by $[0,2]$ is trivial. 
Since (the associated bicharacter of) $q^{-2}_d$ is non-degenerate the subcategory $\V$ is non-degenerate as a braided category. 
Hence by M\"uger's centraliser theorem \cite[Prop.\,4.1]{mu} $\C(N{=}2,d)_{NS}\simeq \T\boxtimes\V$ as ribbon fusion categories.

Finally, we will show that as a tensor category and for odd $d$, $\C(\bZ_d,q^{-2}_{d})$ is equivalent to the category $\V(\bZ_d)$ of $\bZ_d$-graded vector spaces with the trivial associator. The quadratic form $q^{-2}_{d} \in Q(\bZ_d,\bC^*)$ determines the braided tensor structure on $\C(\bZ_d,q^{-2}_{d})$ via the canonical isomorphism from $Q(\bZ_d,\bC^*)$ to the third abelian group cohomology $H^3_{ab}(\bZ_d,\bC^*)$ \cite{js}. The associator on $\C(\bZ_d,q^{-2}_{d})$, i.e.\ the structure as a tensor category, is determined by the image under the homomorphism $H^3_{ab}(\bZ_d,\bC^*)\to H^3(\bZ_d,\bC^*)$. For $d$ odd, this homomorphism is trivial, hence the associator on $\C(\bZ_d,q^{-2}_{d})$ is trivial.

The above discussion is summarised in the following statement.

\bpr\lb{ddn}
For an odd $d$ there is an equivalence of 
	braided
fusion categories $$\C(N{=}2,d)_{NS}\ \simeq\ \T\boxtimes\V(\bZ_d)\ .$$
\epr

\subsection{Universal properties}\label{sec:univ}

Here we formulate universal properties of Temperley-Lieb and pointed fusion categories.
We say that a tensor category $\C$ is {\em freely generated by} an object $X\in\C$ together with a collection of morphisms $\{ f_j:X^{\ot n_j}\to X^{\ot m_j} \}$ making a collection of diagrams $D_s$ commutative if for any tensor category $\D$ the functor of taking values 
$$\Funct_\ot(\C,\D)\to \D',\qquad F\mapsto F(X)$$
is an equivalence. Here, $\Funct_\ot(\C,\D)$ is the category of tensor functors (with tensor natural transformations as morphisms). The target $\D'$ is the category with objects $(Y,\{g_j\})$, where $Y\in\D$ and the $g_j:Y^{\ot n_j}\to Y^{\ot m_j}$ make the collection of diagrams $D_s$, with $X$ replaced by $Y$
and $f_j$ by $g_j$, 
commutative in $\D$. Morphisms $(Y,\{g_j\})\to (Y',\{g'_j\})$ in $\D'$ are morphisms $Y\to Y'$ in $\D$ fitting into commutative squares with all $g_j,g'_j$.
 
\subsubsection{Temperley-Lieb categories}\lb{uptl}

We call an object $T$ of a tensor category $\C$ {\em self-dual} if it comes equipped with morphisms
$$
	n:I\to T \otimes T~~ ,\qquad u:T \otimes T \to I \ ,
$$ 
such that the diagrams 
$$\xymatrix{T \ar[rrr]^{\id} \ar[d]_{\lambda^{-1}_T} &&& T\\ I\ot T\ar[r]^{n\ot \id} & (T\ot T)\ot T \ar[r]^{a_{T,T,T}^{-1}} & T\ot(T\ot T)\ar[r]^{\id\ot u} & T\ot I \ar[u]_{\rho_T}}$$
\be\label{eq:zigzag}
\xymatrix{T \ar[rrr]^{\id} \ar[d]_{\rho^{-1}_T} &&& T\\ T\ot I\ar[r]^{\id\ot n} & T\ot (T\ot T) \ar[r]^{a_{T,T,T}} & (T\ot T)\ot T\ar[r]^{u\ot\id} & I\ot T \ar[u]_{\lambda_T}} 
\ee
commute. 
If there is a scalar $\kappa\in \kk$ such that $u \circ n = \kappa\, \id_I$, we say $T$ has (self-dual) {\em dimension} $\kappa$. 

The category $\TL_\kappa$ freely generated by a self-dual object of non-zero
dimension $\kappa$ is called the {\em Temperley-Lieb category} (see \cite[Chapter XII]{tu}). It has a geometric description as a category with morphism being (isotopy classes of) plane tangles modulo some simple relations. In particular, according to this description the endomorphism algebras $\TL_\kappa(T^{\ot n},T^{\ot n})$ are Temperley-Lieb algebras $TL_n(\kappa)$, i.e.\ algebras with generators $e_i,\ i=1...,n-1$ and relations 
$$e_i^2 = \kappa\ e_i,\qquad e_ie_{i\pm1}e_i = e_i,\qquad e_ie_j = e_je_i\quad |i-j|>1\ .$$

Now write $\kappa=\eta+\eta^{-1}$ for $\eta\in \kk$. 
Consider the simple objects $T_n\in \TL_\kappa$ defined as images of certain idempotents $p_n\in TL_n(\kappa)$ (the {\em Jones-Wenzl projectors}), which are given by, 
	for $n \ge 1$,
$$p_{n+1} = p_n\ot \id - \frac{[n]_\eta}{[n+1]_\eta}(p_n\ot \id) \circ e_n \circ (p_n\ot \id) 
	\quad , \quad p_1 = \id 
\ ,$$
where $[n]_\eta = \frac{\eta^n-\eta^{-n}}{\eta-\eta^{-1}}$ are quantum numbers.
The dimension of $T_n$ (which can be computed as the trace of $p_n$) is $\dim(T_n) = [n{+}1]_\eta$. 
	We set $T_0=I$, the monoidal unit, and from the above definition $T_1=T$ is the generating object.
It is straightforward to see that the endomorphism algebras $\TL_\kappa(T\ot T_n,T\ot T_n)$ are 2-dimensional
	for all $n \ge 1$.

For $\eta$ a root of unity of order $> 2$, the last well-defined Jones-Wenzl projector is $p_{d-1}$, where $d$ is the order of $\eta$ if it is odd and half the order of $\eta$ if it is even.
In this case the category $\TL_\kappa$ has a maximal fusion quotient $\T_\kappa$ which can be defined as the quotient 
$$
\T_\kappa ~:=~ \TL_\kappa/\langle p_{d-1}\rangle
$$ 
by the ideal of morphisms tensor generated by the Jones-Wenzl projector $p_{d-1}\in TL_{d-1}(\kappa)$, see \cite{eo}. 
Moreover the ideal of morphisms tensor generated by the Jones-Wenzl projector $p_{d-1}$ is the unique 
non-zero proper
tensor ideal in $\TL_\kappa$ \cite{gw}, that is, any non-faithful tensor functor $\TL_\kappa\to\D$ factors through $\T_\kappa\to\D$. 
Thus we have the following.
\bth
A tensor functor from $\T_\kappa$ to a tensor category $\D$ is determined by a self-dual object of dimension $\kappa$ in $\D$ with vanishing Jones-Wenzl projector $p_{d-1}$. 
\eth

The next corollary provides an easy-to-use replacement for the vanishing
condition on the Jones-Wenzl projector.
Recall that the simple objects of $\T_\kappa$ are $T_i,\ i=0,...,d-2$
with $T_0=I, T_1=T$. The tensor product with $T$ is $T\ot T_i\simeq
T_{i-1}\op T_{i+1}$ for $0<i<d-2$ and $T\ot T_{d-2}\simeq T_{d-3}$.

\bco\lb{ftl}
Let $\D$ be a rigid fusion category with simple objects $S_i,\
i=0,...,d-2$ and the tensor product  $S_1\ot S_i\simeq S_{i-1}\op
S_{i+1}$ for $0<i<d-2$ and $S_1\ot S_{d-2}\simeq S_{d-3}$. A tensor
functor $\TL_\kappa \to \D$ such that $T_i \mapsto S_i$ factors through
$\T_\kappa$.
\eco

\bpf
The non-faithfulness of the tensor functor is manifest since
$\TL_\kappa(T\ot T_{d-2},T\ot T_{d-2})$ is 2-dimensional, while
$\D(S_1\ot S_{d-2},S_1\ot S_{d-2})$ is only 1-dimensional.
\epf

See also \cite{da} for details.

\subsubsection{The categories $\C(N{=}2,d)_{NS}$ for odd $d$}

Here, we describe a universal property of $\C(N{=}2,d)_{NS}$ for odd $d$ as a tensor category.
	This description makes use of group actions on tensor categories and equivariant objects, which we review in Appendix \ref{sec:equiv+pointed}. In the following proposition, a pointed subcategory of a tensor category $\D$ with underlying group $\bZ_d$ acts by conjugation, and  $\D^{\bZ_d}$ denotes the corresponding tensor category of equivariant objects.

\bth\lb{upn}
Let $d$ be odd. A tensor functor $F:\C(N{=}2,d)_{NS}\to\D$ is determined by 
\begin{itemize}
\item a tensor functor $\V(\bZ_d)\to \D$, 
\item a self-dual object $T=F([1,d])$ in the category $\D^{\bZ_d}$ of quantum dimension $\dim(T)=2\cos\big(\frac{\pi}{d}\big)$ such that the induced functor $\TL_{2\cos(\frac{\pi}{d})}\to \D^{\bZ_d}$ is not faithful.
\end{itemize}
\eth
\bpf
By Proposition \ref{ddn}, the category $\C(N{=}2,d)_{NS}$ is tensor equivalent to the Deligne product $\T\boxtimes\V(\bZ_d).$
By Theorem \ref{ciz}, a tensor functor $F:\T\boxtimes\V(\bZ_d)\to\D$ is determined by a tensor functor $\V(\bZ_d)\to \D$ and a tensor functor $\T\to \D^{\bZ_d}$.

The dimension of $[1,d]\in\T$ (which coincides with the dimension of $[1,0,0]$ in $\E$)
is equal to $2\cos\big(\frac{\pi}{d}\big)$.
The fusion rules of $\T$ (see \eqref{eq:E-fusion}) show that it is freely generated as a tensor category by $[1,d]$, and that the Jones-Wenzl projector $p_{d-1}$ vanishes (Corollary \ref{ftl}). By semi-simplicity, it follows that $\TL_{2\cos(\frac{\pi}{d})}\to \T$ descends to a tensor equivalence $\T_{2\cos(\frac{\pi}{d})}\to \T$. Consequently,
a tensor functor $\T\to \D^{\bZ_d}$ is determined by a self-dual object $T=F([1,d])$ in the category $\D^{\bZ_d}$ with quantum dimension $\dim(T)=2\cos\big(\frac{\pi}{d}\big)$ and such that the induced functor $\TL_{2\cos(\frac{\pi}{d})}\to \D^{\bZ_d}$ is not faithful.
\epf

\section{Matrix factorisations}\label{sec:mf}

\subsection{Categories of matrix factorisations and tensor products}\label{sec:cat-mf}

A {\em matrix factorisation} over a commutative $\kk$-algebra $S$ of an element
$W \in S$ is a $\bZ_2$-graded free $S$-module $M$ together with a {\em twisted differential} $d^M : M \to M$ of odd degree satisfying $d^M\circ d^M = W$. Here, the right hand side stands for the endomorphism $m \mapsto W.m$. 
We will often omit the superscript $M$ in $d^M$ and display the $\bZ_2$-grading explicitly as $M=M_0 \oplus M_1$, $d = d_0 \oplus d_1$ or graphically as 
$$
M~:\quad
\xymatrix{M_1 \ar@/^10pt/[rr]^{d_1} && M_0 \ar@/^10pt/[ll]^{d_0}}\quad .
$$
A matrix factorisation is of {\em finite rank} if its underlying free $S$-module is of finite rank. 

We distinguish several categories of matrix factorisations:
\begin{itemize}
\item 
$\MF_{S,W}$: 
Objects are matrix factorisations $M = (M,d)$ and the morphism space $\MF_{S,W}(M,N)$ consists of all $S$-linear maps from $M$ to $N$. The $\bZ_2$-grading of $M$ and $N$ induces a $\bZ_2$-grading on $\MF_{S,W}(M,N)$.
The twisted differentials of $M$ and $N$ combine to a (non-twisted) degree 1 differential $\delta$ on $\MF(M,N)$ given by $\delta(f) = d^N \circ f - (-1)^{|f|} f \circ d^M$, where $|f|$ is the $\bZ_2$-degree of $f$. In this way,
the morphisms in $\MF_{S,W}$ form a $\bZ_2$-graded complex.
\item 
$\ZMF_{S,W}$: 
Objects are as for $\MF_{S,W}$ and morphisms from $M$ to $N$ are degree zero cycles in $\MF_{S,W}(M,N)$, that is
$$
	\ZMF_{S,W}(M,N) = \{ f : M \to N | \text{ $f$ is $S$-linear of degree 0 and $\delta(f)=0$ } \} \ .
$$
\item 
$\HMF_{S,W}$: 
Objects are as for $\MF_{S,W}$ and the set of morphisms from $M$ to $N$ is the degree zero homology in $\MF_{S,W}(M,N)$, that is
$$
	\HMF_{S,W}(M,N) = \ZMF_{S,W}(M,N) \, / \, \{  \delta(g) \, | \text{ $g : M \to N$ is $S$-linear of degree 1 } \} \ .
$$
\end{itemize}
We will often write morphisms $f \in \ZMF_{S,W}(M,N)$ (or representatives of classes in $\HMF_{S,W}(M,N)$) in a diagram as follows:
$$
  \xymatrix{M_1 \ar@/^10pt/[rr]^{d_1^M} \ar[dd]_{f_1} && M_0 \ar@/^10pt/[ll]^{d_0^M}  \ar[dd]^{f_0} \\ \\
N_1 \ar@/^10pt/[rr]^{d_1^{N}} && N_0 \ar@/^10pt/[ll]^{d_0^{N}}
}
$$
That $f$ is in $\ZMF(M,N)$ is equivalent to $f_0$ and $f_1$ being $S$-linear maps such that the subdiagram with upward curved arrows commutes and that with downward curved arrows commutes: 
$$ 
f_0 \circ d_1^M = d_1^N \circ f_1 \quad , \quad 
f_1 \circ d_0^M = d_0^N \circ f_0 \ .
$$
In fact, if $W$ is not a zero-divisor in $S$, one condition implies the other (see \cite[Ch.\,7]{yoshinobuch}).

For more on matrix factorisations in general we refer to foundational works \cite{eisenbud,buchweitz} or for example to \cite{yoshinobuch,khovroz}.

\medskip

The above definitions can be made also for bimodules, giving rise to the notion of a matrix bifactorisation \cite{cr1}.

\begin{defin}
A {\em matrix bifactorisation} over $S$ of
$W$ is a pair $\left( M,d^M \right)$ where $M$ is a $\bZ_2$-graded free $S$-$S$-bimodule and $d^M:M \rightarrow M$ an $S$-$S$-bimodule endomorphism of degree 1 satisfying $d^M \circ d^M=W.\id_M-\id_M.W$, where the right hand side stands for the map $m \mapsto W.m - m.W$.
\end{defin}

Here, an $S$-$S$-bimodule is called free if it is free as an $S \otimes_\kk S$-left module. 
As with matrix factorisations one can define morphisms of matrix bifactorisation (in this case morphisms of bimodules instead of simply modules). We denote the resulting differential $\bZ_2$-graded category as $\MF_{\mathrm{bi};S,W}$. The associated categories with morphisms which are degree zero cycles and degree zero homology classes are defined as before and will be denoted as $\ZMF_{\mathrm{bi};S,W}$ and $\HMF_{\mathrm{bi};S,W}$.

As the algebra $S$ and the element
$W$ will be clear from the context (in fact, we will soon restrict to $S = \bC[x]$ and $W = x^d$), we will omit the subscript $S,W$ from now on.

\medskip

For $S=\kk[x_1,\dots,x_n]$ and $W \in S$ a {\em potential} (i.e.\ $Jac(W)$ is finite dimensional, see \cite{khovroz} for details),
the category $\HMFbi$ is tensor
\cite{cr1,cm1} (for $S$ arbitrary, it is still non-unital tensor).
The tensor product of $M,N \in \MFbi$ is given by 
$$
	M \otimes_S N 
	\quad , \quad d = d^M \otimes_S \id_N + \id_M \otimes_S d^N \ .
$$
In the following we will just write $\otimes$ for $\otimes_S$.
The above definition hides a Koszul sign: for $m \in M$ and $n \in N$ we have $(\id_M \otimes d^N)(m \otimes n) = (-1)^{|m|} \, m \otimes d^N(n)$, where $|m| \in  \bZ_2$ denotes the degree of $m$. Thus, if we spell out the twisted differential of $M \otimes N$ in components and make the Koszul sign explicit, we have
$$
M \otimes N~:\quad
\xymatrix{
{\begin{array}{c} {M_1 \otimes N_0} \\  {\oplus} \\ {M_0 \otimes N_1} \end{array} }
\ar@/^10pt/[rrrrrr]^{d_1^{M\otimes N}=\left( \begin{matrix} d_1^M \otimes \id_{N_0} & \id_{M_0} \otimes d_1^{N} \\ -\id_{M_1} \otimes d_0^{N} & d_0^M \otimes \id_{N_1} \end{matrix} \right)} &&&&&& 
{\begin{array}{c} M_0 \otimes N_0 \\ \oplus  \\ M_1 \otimes N_1 \end{array} } 
\ar@/^10pt/[llllll]^{d_0^{M\otimes N}=\left( \begin{matrix} d_0^M \otimes \id_{N_0} & -\id_{M_1} \otimes d_1^{N} \\ \id_{M_0} \otimes d_0^{N} & d_1^M \otimes \id_{N_1} \end{matrix}\right)}}
$$ 
The associativity isomorphisms are simply those of the underlying tensor category of bimodules. However, the unit object in the category of $\bZ_2$-graded $S$-$S$-bimodules, the bimodule $S$, is not free as an $S \otimes_\kk S$-left module. As a consequence, the categories $\MFbi$ and $\ZMFbi$ are non-unital tensor. On the other hand, $\HMFbi$ has a unit object, which we give explicitly in the case $S=\bC[x]$ and $W=x^d$ below. For the general case we refer to 
\cite{cr1,cm1}.
For more on tensor products see \cite{yoshino,khovroz,brunnerrogg,cr1,Dyckerhoff:2011vf,Carqueville:2011zea,cm1,murfet}. 

\medskip

{}From here on and for the remainder Section \ref{sec:mf} we fix
$$
	S = \bC[x] \quad , \qquad W=x^d \quad , \quad \text{where} \quad
	d \in \bZ \quad , \quad d \ge 2 \ .
$$
For calculations it will often be convenient to describe $\bC[x]$-$\bC[x]$-bimodules as $\bC[x,y]$-left modules $M$. Here, the left action of $p \in \bC[x]$ is by acting on $M$ with $p(x)$ and the right action by acting with $p(y)$. We will employ this tool without further mention.

The tensor unit in $\HMFbi$ is
$$
I ~:\quad\xymatrix{\bC[x,y] \ar@/^10pt/[rrr]^{d_1 = x-y} &&& \bC[x,y] \ar@/^10pt/[lll]^{d_0 = \frac{x^d-y^d}{x-y}}} \quad .
$$
The left and right unit isomorphisms $\lambda_M : I \otimes M \to M$ and $\rho_M : M \otimes I \to M$ are given by
$$\xymatrix{
I\ot M \ar[dd]_{\lambda_M} && 
{\begin{array}{c} {I_1 \otimes M_0} \\  {\oplus} \\ {I_0 \otimes M_1} \end{array} }
\ar[dd]_{(0\ L_{M_1})} \ar@/^10pt/[rrrr]^{\left(\begin{array}{rr} \scriptstyle{(x-y)\ot\id} & {\scriptstyle{1\ot d_1^M}} \\ {\scriptstyle{-1\ot d_0^M}} &  \scriptstyle{\frac{x^d-y^d}{x-y}\ot\id} \end{array}\right)} &&&& 
{\begin{array}{c} {I_0 \otimes M_0} \\  {\oplus} \\ {I_1 \otimes M_1} \end{array} }
\ar[dd]^{(L_{M_0}\ 0)} \ar@/^10pt/[llll]^{\left(\begin{array}{rr}  \scriptstyle{\frac{x^d-y^d}{x-y}\ot\id} & {\scriptstyle{-1\ot d_1^M}} \\ {\scriptstyle{1\ot d_0^M}} & \scriptstyle{(x-y)\ot\id} \end{array}\right)}
 \\ \\ 
M  & & M_1 \ar@/^10pt/[rrrr]^{d_1^M}  &&&& M_0 \ar@/^10pt/[llll]^{d_0^M} 
}$$
\be\label{eq:unit-isos}
\xymatrix{
M\ot I \ar[dd]_{\rho_M} && 
{\begin{array}{c} {M_1 \otimes I_0} \\  {\oplus} \\ {M_0 \otimes I_1} \end{array} }
\ar[dd]_{(R_{M_1}\ 0)} \ar@/^10pt/[rrrr]^{\left(\begin{array}{rr} {\scriptstyle{d_1^M}\ot 1} & \scriptstyle{\id\ot(x-y)}  \\  \scriptstyle{-\id\ot\frac{x^d-y^d}{x-y}} & {\scriptstyle{d_0^M\ot 1}}  \end{array}\right)} &&&& 
{\begin{array}{c} {M_0 \otimes I_0} \\  {\oplus} \\ {M_1 \otimes I_1} \end{array} }
\ar[dd]^{(R_{M_0}\ 0)} \ar@/^10pt/[llll]^{\left(\begin{array}{rr}  {\scriptstyle{d_0^M\ot 1}} & \scriptstyle{-\id\ot(x-y)} \\  \scriptstyle{\id\ot\frac{x^d-y^d}{x-y}} & {\scriptstyle{d_1^M\ot 1}} \end{array}\right)}
 \\ \\
M  & & M_1 \ar@/^10pt/[rrrr]^{d_1^M}  &&&& M_0 \ar@/^10pt/[llll]^{d_0^M} 
}\ee	
The maps $L$ and $R$ are, for a given $\bC[x]$-$\bC[x]$-bimodule $N$, defined as
\begin{align*}
	L_N : \bC[x,y] \otimes N &\longrightarrow N 
	&
	R_N : N \otimes \bC[x,y] &\longrightarrow N 
	\\
	f(x,y) \otimes n &\longmapsto f(x,x).n
	&
	n \otimes f(x,y) &\longmapsto n.f(x,x)
	\nonumber
\end{align*}
It is easy to verify that $\lambda_M$ and $\rho_M$ are in $\ZMFbi$. With some more work, one sees that they have homotopy inverses, see \cite{cr1}.

\medskip

Finite rank factorisations in $\HMFbi$ have right duals \cite{cr2,cm1}.
We will only need explicit duals of matrix factorisations $M \in \HMFbi$ for which $M_0$ and $M_1$ are of rank 1. In this case we have \cite{cr2}
$$
M~:\quad
\xymatrix{\bC[x,y] \ar@/^10pt/[rr]^{d_1(x,y)} && \bC[x,y] \ar@/^10pt/[ll]^{d_0(x,y)}}
\qquad \leadsto \qquad
M^+~:\quad
\xymatrix{\bC[x,y] \ar@/^10pt/[rr]^{d_1^{M^+} := -d_1(y,x)} &&\bC[x,y] \ar@/^10pt/[ll]^{d_0^{M^+} := d_0(y,x)}}\quad .
$$
Note that $I^+ = I$. Since the corresponding duality maps play an important role in our construction, we take some time to recall their explicit form and some properties from \cite{cr2}. The coevaluation  $coev_M : I \to M \otimes M^+$ is the simpler of the two,
$$
\xymatrix{
I \ar[ddd] & & \bC[x,z] \ar@/^10pt/[rrrr]^{x-z} \ar[ddd]_{\left(\begin{array}{c} \scriptstyle{1} \\ \scriptstyle{1} \end{array}\right)} &&&& \bC[x,z] \ar@/^10pt/[llll]^{\frac{x^d-z^d}{x-z}} \ar[ddd]^{\left(\begin{array}{r} \frac{d_1(x,y)-d_1(z,y)}{x-z} \\ \frac{d_0(x,y)-d_0(z,y)}{x-z} \end{array}\right)} \\ \\ \\
M \otimes M^+ && \bC[x,y,z]^{\op 2} \ar@/^10pt/[rrrr]^{\left(\begin{array}{rr} \scriptstyle{d_1(x,y)} & \scriptstyle{-d_1(z,y)} \\ \scriptstyle{-d_0(z,y)} &  \scriptstyle{d_0(x,y)} \end{array}\right)} &&&& \bC[x,y,z]^{\op 2} \ar@/^10pt/[llll]^{\left(\begin{array}{rr} \scriptstyle{d_0(x,y)} & \scriptstyle{d_1(z,y)} \\ \scriptstyle{d_0(z,y)} & \scriptstyle{d_1(x,y)} \end{array}\right)}
}
$$
Here the left and the right bottom instances of $\bC[x,y,z]^{\op 2}$ correspond to 
$$(M \otimes M^+)_1 = {\begin{array}{c} {M_1 \otimes M^+_0} \\  {\oplus} \\ {M_0 \otimes M^+_1} \end{array} }\qquad ,
\qquad (M \otimes M^+)_0 = {\begin{array}{c} {M_0 \otimes M^+_0} \\  {\oplus} \\ {M_1 \otimes M^+_1} \end{array} }
\quad , 
$$ respectively.
It is immediate that this is indeed a morphism in $\ZMFbi$. The evaluation $ev_M : M^+ \otimes M \to I$ takes the form
$$
\xymatrix{
M^+ \otimes M \ar[ddd] && \bC[x,y,z]^{\op 2} \ar[ddd]_{(B_M\ C_M)} \ar@/^10pt/[rrrr]^{\left(\begin{array}{rr} \scriptstyle{-d_1(y,x)} & \scriptstyle{d_1(y,z)} \\ \scriptstyle{-d_0(y,z)} &  \scriptstyle{d_0(y,x)} \end{array}\right)} &&&& \bC[x,y,z]^{\op 2} \ar[ddd]^{(A_M\ 0)} \ar@/^10pt/[llll]^{\left(\begin{array}{rr} \scriptstyle{d_0(y,x)} & \scriptstyle{-d_1(y,z)} \\ \scriptstyle{d_0(y,z)} &  \scriptstyle{-d_1(y,x)} \end{array}\right)}
 \\ \\ \\
I  & & \bC[x,z] \ar@/^10pt/[rrrr]^{x-z}  &&&& \bC[x,z] \ar@/^10pt/[llll]^{\frac{x^d-z^d}{x-z}} 
}
$$
Here the left and the right top instances of $\bC[x,y,z]^{\op 2}$ correspond to 
$$(M^+ \otimes M)_1 = {\begin{array}{c} {M^+_1 \otimes M_0} \\  {\oplus} \\ {M^+_0 \otimes M_1} \end{array} }\qquad,\qquad (M^+ \otimes M)_0 = {\begin{array}{c} {M^+_0 \otimes M_0} \\  {\oplus} \\ {M^+_1 \otimes M_1} \end{array} }\quad,$$ respectively.
The $\bC[x,z]$-module maps $A_M, B_M, C_M$ are defined as follows. The map $C_M$ is simply minus the projection onto terms independent of $y$: $C_M(y^m) = -\delta_{m,0}$. For $A_M$ and $B_M$ we introduce the auxiliary function
$$
	\mathcal{G}_M(f) = \frac{1}{2 \pi i} \oint \frac{x-z-y}{y \, d_1(y,z)} \, f(x,y,z) dy 
	\qquad , \quad f \in \bC[x,y,z] \ .
$$
The contour integration is along a counter-clockwise circular contour enclosing all poles. It is not immediately evident but still true that $\mathcal{G}_M(f)$ is a polynomial. One way to see this is to rewrite $\mathcal{G}_M(f) = \frac{1}{2 \pi i} \oint \frac{x-z-y}{y (y^d-z^d)} \, d_0(y,z)f(x,y,z) dy$ and to expand $(y^d-z^d)^{-1} = \sum_{m=0}^\infty (z/y)^m$. In this way one can rewrite the integrand as a formal Laurent series in $y$ whose coefficients are polynomials in $x,z$. The contour integration picks out the coefficient of $y^{-1}$. 

We will need two further properties of $\mathcal{G}_M$:
\be\label{eq:G-properties}
	\mathcal{G}_M\big(d_1(y,z) \, y^m\big) = (x-z) \delta_{m,0}
		\quad , \quad
	\mathcal{G}_M\big(d_1(y,x) \, f(x,y,z)\big) \in (x-z) \bC[x,z] \ .
\ee
The first property is clear. For the second property, let $g(x,z) := \mathcal{G}_M\big(d_1(y,x) \, f(x,y,z)\big)$. The condition $g(z,z)=0$ is then immediate from the first property. 

We can now give the maps $A_M$ and $B_M$:
$$
	A_M(f) = - \mathcal{G}_M(f)
	\quad , \quad
	B_M(f) = \frac{\mathcal{G}_M\big(d_1(y,x)f(x,y,z)\big)}{x-z} \ .
$$
To verify that  $ev_M \in \ZMFbi(M^+ \otimes M,I)$, it suffices to check $(ev_M)_0 \circ d^{M^+ \otimes M}_1 = d^I_1 \circ (ev_M)_1$ on $(y^m,y^n)$ for all $m,n \ge 0$. This is straightforward using \eqref{eq:G-properties}:
\begin{align*}
(ev_M)_0 \circ d^{M^+ \otimes M}_1(y^m,y^n) &= A_M\big( -d_1(y,x) y^m + d_1(y,z) y^n \big)
=  \mathcal{G}_M(d_1(y,x)y^m) - (x-z) \delta_{n,0} \ ,
\nonumber\\
d^I_1 \circ (ev_M)_1(y^m,y^n) &= (x-z)(B_M(y^m) + C_M(y^n)) = \mathcal{G}_M(d_1(y,x)y^m) - (x-z) \delta_{n,0} \ .
\nonumber
\end{align*}

The zig-zag identities for $ev_M$ and $coev_M$ are verified in \cite[Thm.\,2.5]{cr2}.

\subsection{Permutation type matrix bifactorisations}\label{sec:permutation-mf}

We fix the primitive $d$'th root of unity\footnote{
Anticipating Remark \ref{rem:galois} below, the reader may check that all statements below -- except for Theorem \ref{thm:main} -- work equally for any other choice of primitive $d$'th root of unity.}
$$
	\eta = e^{\frac{2 \pi i}{d}} \ .
$$ 
For a subset $S \subset \bZ_d$ write $\overline S = \bZ_d \setminus S$. By a {\em permutation type matrix bifactorisations} we mean
\be\label{pfs}
P_S ~:\quad\xymatrix{\bC[x,y] \ar@/^10pt/[rrr]^{d_1 = \prod\limits_{j \in S}(x-\eta^j y)} &&& \bC[x,y] \ar@/^10pt/[lll]^{d_0 = \prod\limits_{j \in \overline S}( x-\eta^j y)}} \quad .
\ee
For example, $I = P_{\{0\}}$. The bifactorisations $P_\emptyset$ and $P_{\{0,1,\dots,d-1\}}$ are isomorphic to the zero object in $\HMFbi$. The remaining $P_S$ are non-zero and mutually distinct.
To see this, in the following remark we recall a useful tool from $\cite{khovroz}$.

\bre\label{rem:H(M)}
Given a matrix bifactorisation $(M,d)$, we obtain a $\bZ_2$-graded complex by considering the differential $\bar d$ on $M / \langle x,y \rangle M$. Since $x^d-y^d \in \langle x,y \rangle$, $\bar d$ is indeed a differential. Denote by $H(M)$ the homology of this complex. Then \cite[Prop.\,8]{khovroz} states that $f \in \HMFbi(M,N)$ is an isomorphism in $\HMFbi$ if and only if the induced map $H(f) : H(M) \to H(N)$ is an isomorphism of $\bC$-vector spaces.
\ere

\ble
Let $R,S \subset \bZ_d$ be nonempty proper subsets. The permutation type matrix bifactorisations $P_R$ and $P_S$ are non-zero, and they are isomorphic in $\HMFbi$ if and only if $R=S$.
\ele

\bpf
For a non-empty proper subset $S$, the matrix factorisation $P_S$ is reduced, that is, the differential $\bar d$ induced on the quotient $P_S/\langle x,y\rangle P_S$ is zero. Thus $H(P_S) \simeq \bC \oplus \bC$. It follows that $f \in \ZMFbi(P_S,P_R)$ is an isomorphism in $\HMFbi$ if and only if $f_0$ and $f_1$ contain a non-zero constant term. Writing out the condition that $f$ is a cycle shows that this is possible only for $R=S$.
\epf

We will mostly be concerned with a special subset of permutation type bifactorisations, namely those with consecutive index sets. For $a \in \bZ_d$ and $\lambda \in \{0,1,2,\dots,d-2\}$ we write
$$
	P_{a:\lambda} := P_{\{a,a+1,\dots,a+\lambda\}} \ .
$$
We define $\Pd$ to be the full subcategory of $\HMFbi$ consisting of objects isomorphic (in $\HMFbi)$ to finite direct sums of the $P_{a:\lambda}$. A key input in our construction is the following result established in \cite[Sect.\,6.1]{brunnerrogg}.

\bth\label{thm:br-decomp}
$\Pd$ is closed under taking tensor products. Explicitly, for $\lambda,\mu \in \{0,\dots,d-2\}$,
$$
	P_{m:\lambda}  \otimes P_{n:\mu} \simeq \bigoplus_{\nu=|\lambda-\mu|~\mathrm{step}\,2}^{\min(\lambda+\mu,2d-4-\lambda-\mu)} P_{m+n-\frac12(\lambda+\mu-\nu):\nu} \quad .
$$
\eth

For the dual of a permutation type matrix bifactorisations one finds $(P_S)^+ \simeq P_{-S}$. Explicitly:
\be\label{eq:PS+P-S-iso}
\xymatrix{
 P_{-S} \ar[dd] &&&& \bC[x,y] \ar@/^10pt/[rrr]^{\prod\limits_{j \in S}(x-\eta^{-j} y)} \ar[dd]_{(-1)^{|S|+1}\prod_{j\in S}\eta^{-j}} &&&\bC[x,y] \ar@/^10pt/[lll]^{\prod\limits_{j \in \overline S}(x-\eta^{-j} y)}  \ar[dd]^{1} \\ \\
(P_S)^+ &&&& \bC[x,y] \ar@/^10pt/[rrr]^{-\prod\limits_{j \in S}(y-\eta^{j}x)} &&& \bC[x,y] \ar@/^10pt/[lll]^{\prod\limits_{j \in \overline S}(y-\eta^{j}x)}
}
\ee
The self-dual permutation type matrix bifactorisations of the form $P_{a:\lambda}$ therefore have to satisfy $2a \equiv -\lambda \mod d$. Depending on the parity of $d$, one finds:
\begin{itemize}
\item $d$ even: $\lambda$ must be even and $a \equiv \frac\lambda2 \mod d$ or $a \equiv \frac{\lambda+d}2 \mod d$,
\item $d$ odd: $\lambda$ can be arbitrary and $a \equiv \frac{d-1}2 \lambda \mod d$.
\end{itemize}

\subsection{A tensor functor from $\bZ_d$ to $\Pd$}\label{sec:Zd_to_Pd}

Consider the algebra automorphism $\sigma$ of $\bC[x]$ which acts on $x$ as $\sigma(x) = \eta x$. It leaves the potential $x^d$ invariant and generates the group of algebra automorphisms with this property. We get a group isomorphism
$$
	\bZ_d ~\longrightarrow~ \mathrm{Aut}(\bC[x] \text{ with $x^d$ fixed})
	\quad , \quad k \mapsto \sigma^k \ .
$$
Given a matrix bifactorisation $M \in \MFbi$ and $a,b \in \bZ_d$, we denote by ${}_aM_b$ the matrix bifactorisation whose underlying $\bC[x]$-bimodule is equal to $M$ as a $\bZ_2$-graded $\bC$-vector space, but has twisted left/right actions ($p \in \bC[x]$, $m \in M$):
$$
	(p,m) \mapsto \sigma^{-a}(p).m \quad , \quad
	(m,p) \mapsto m.\sigma^{b}(p) \ ,
$$
where the dots denotes the left/right action on the original bimodule $M$. Since $\bZ_d$ is abelian, we get a left action even if we were to omit the minus sign in $\sigma^{-a}$, but we include it to match the conventions of \cite[Sect.\,7.1]{cr3}.
For permutation type matrix bifactorisations we have isomorphisms:
\be\label{eq:PS-with-twisted-action-iso-to-PS}
  \xymatrix{
P_{S-a-b}  \ar[dd]_{s_{a,b}} && \bC[x,y] \ar@/^10pt/[rrr]^{\prod\limits_{j \in S}(x-\eta^{j-a-b} y)}  \ar[dd]_{\eta^{-|S|a} \,\cdot\, \sigma^{-a} \otimes \sigma^{b} } &&& \bC[x,y] \ar@/^10pt/[lll]^{\prod\limits_{j \in \overline S}(x-\eta^{j-a-b} y)}  \ar[dd]^{ \sigma^{-a} \otimes \sigma^{b}}
  \\ \\
 {}_a(P_S)_b && {}_a(\bC[x,y])_b \ar@/^10pt/[rrr]^{\prod\limits_{j \in S}(x-\eta^j y)} &&& {}_a(\bC[x,y])_b \ar@/^10pt/[lll]^{\prod\limits_{j \in \overline S}(x-\eta^j y)}  
}
\ee
Here, $\sigma^{-a} \otimes \sigma^{b}$ is the automorphism of $\bC[x,y]$ which acts as $x \mapsto \eta^{-a} x$ and $y \mapsto \eta^{b}y$. 

The following lemma is straightforward.

\ble\label{lem:a()b-functor}
For all $a,b \in \bZ_d$, 
${}_a(-)_{b}$ defines an auto-equivalence of $\HMFbi$ and of $\Pd$. If $b=-a$, this auto-equivalence is 
tensor with ${}_a(M \otimes N)_{-a} = {}_aM_{-a} \otimes {}_aN_{-a}$ and $s_{a,-a} : I \to {}_aI_{-a}$.
\ele

Consider the objects ${}_aI \in \HMFbi$ for $a \in \bZ_d$. Applying the functor ${}_a(-)$ to the unit isomorphism $\lambda_{{}_bI} : I \otimes {}_bI \to {}_bI$ gives the isomorphism
\be\label{eq:mu-morph}
	\mu_{a,b}  := {}_a(\lambda_{{}_bI}) ~:~ {}_aI \otimes {}_bI \to 
	{}_{a+b}I \ .
\ee
By $\underline{\bZ_d}$ we mean the monoidal category whose set of objects is $\bZ_d$, whose set of morphisms consists only of the identity morphisms, and whose tensor product functor is the group operation (i.e.\ addition), see Appendix \ref{ceo}. 

\bpr\label{prop:chi-functor}
$\chi : \underline{\bZ_d} \to \Pd$, $\chi(a) = {}_aI$, together with 
$\mu_{a,b} : \chi(a) \otimes \chi(b) \to \chi(a+b)$, 
defines a tensor functor.
\epr

\bpf
First note that by \eqref{eq:PS-with-twisted-action-iso-to-PS}, ${}_aI \simeq P_{\{-a\}}$, so that indeed $\chi(a) \in \Pd$.
It is shown in \cite[Prop.\,7.1]{cr3} that the $\mu_{a,b}$ satisfy the associativity condition
$$
	\mu_{a,b+c} \circ (\id_{{}_aI} \otimes \mu_{b,c}) = \mu_{a+b,c} \circ (\mu_{a,b} \otimes \id_{{}_cI}) 
	\qquad \text{for all} \quad a,b,c \in \bZ_d \ .
$$
This amounts to the hexagon condition for the coherence isomorphisms $\mu_{a,b}$.
\epf

We can now construct two tensor functors $\underline{\bZ_d} \to \mathrm{Aut}_\otimes(\Pd)$. The first functor takes $a \in \bZ_d$ to ${}_a(-)_{-a}$; we denote this functor by $A$. This functor is strictly tensor: $A(0)=\Id$ and $A(a) \circ A(b) = A(a+b)$.

The second functor is the adjoint action of $\chi$; we denote it by $\Ad_\chi$. Given $a\in\bZ_d$, on objects the functor $\Ad_\chi(a)$ acts as $M \mapsto \chi(a) \otimes M \otimes \chi(-a)$. Morphism $f : M \to N$ get mapped to $\id_{\chi(a)} \otimes f \otimes \id_{\chi(-a)}$. The isomorphisms $\mu_{-a,a} : \chi(-a) \otimes \chi(a) \to \chi(0) = I$ give the tensor structure on $\Ad_\chi(a)$. 
So far we saw that for all $a \in \bZ_d$, $\Ad_\chi(a) \in  \mathrm{Aut}_\otimes(\Pd)$. Next we need the coherence isomorphisms $\Ad_\chi(a) \circ \Ad_\chi(b) \to \Ad_\chi(a+b)$. These are simply given by $\mu_{a,b} \otimes (-) \otimes \mu_{-b,-a}$. 

The following lemma will simplify the construction of $\bZ_d$-equivariant structures below.

\ble\label{lem:A-iso-Ad}
$A$ and $\Ad_\chi$ are naturally isomorphic as tensor functors.
\ele

\bpf
We need to provide a natural monoidal isomorphism $\alpha : \Ad_\chi \to A$. That is, for each $a \in \bZ_d$ we need to give a natural monoidal isomorphism $\alpha(a) : \Ad_\chi(a) \to A(a)$, such that the diagram
\be\label{eq:alpha-mon}
	\xymatrix{ 
	\Ad_\chi(a) \circ \Ad_\chi(b) \ar[d]_-{\alpha(a) \circ \alpha(b)} \ar[r]^-{\mu_*} & \Ad_\chi(a+b) \ar[d]^-{\alpha(a+b)}
	\\
	 A(a) \circ A(b)  \ar@{=}[r] & A(a+b)
	}
\ee
commutes, where $\mu_* := \mu_{a,b} \otimes (-) \otimes \mu_{-b,-a}$. Define
$$
	\alpha(a)_M 
	:=
	\big[ {}_aI \otimes M \otimes {}_{-a}I 
		\xrightarrow{ {}_a(\lambda_M) \otimes (s^{-1}_{-a,a})_{-a} }
		{}_aM \otimes I_{-a}
		\xrightarrow{ {}_a( \rho_M )_{-a} }
		{}_aM_{-a}
	\big] \ .
$$

\noindent
{\em $\alpha(a)$ is tensor:} We need to verify commutativity of
$$
	\xymatrix{ 
	\Ad_\chi(a)(M) \otimes \Ad_\chi(a)(N) \ar[d]_-{\alpha(a)_M \otimes \alpha(a)_N} \ar[r]^-{\sim} & \Ad_\chi(a)(M \otimes N) \ar[d]^-{\alpha(a)_{M \otimes N}}
	\\
	 A(a)(M) \otimes A(a)(N)  \ar@{=}[r] & A(a)(M \otimes N)
	}
$$
where the top isomorphism is $\id \otimes \mu_{-a,a} \otimes \id$. Commutativity of this diagram is a straightforward calculation if one notes the following facts: $M \otimes {}_{-a}I = M_{-a} \otimes I$ and $M_{-a} \otimes {}_aN = M \otimes N$ (equal as matrix factorisations, not just isomorphic), and
$$
	\big[ M \otimes {}_{-a} I \xrightarrow{\id \otimes (s^{-1}_{-a,a})_{-a}} M \otimes I_{-a} \xrightarrow{ (\rho_M)_a} M_{-a} \big] ~=~ \big[ M_{-a} \otimes I \xrightarrow{ \rho_{M_{-a}} } \big] \ .
$$

\noindent
{\em $\alpha$ satisfies \eqref{eq:alpha-mon}:}
One way to see this is to act on elements. The unit isomorphisms \eqref{eq:unit-isos} are non-zero only on summands in the tensor products involving $I_0$, in which case they act as
$$
	\lambda_M ~:~ p(x,y) \otimes m \mapsto p(x,x).m
	\quad , \quad
	\rho_M ~:~ m \otimes p(x,y)\mapsto m.p(x,x) \ .
$$
One verifies that the top and bottom path in \eqref{eq:alpha-mon} amount to mapping
$$
	u(x,y) \otimes v(x,y) \otimes m \otimes p(x,y) \otimes q(x,y)
	~\in~ ({}_aI)_0 \otimes ({}_bI)_0 \otimes M \otimes ({}_{-b}I)_0 \otimes ({}_{-a}I)_0
$$
to
$$
	\big\{ \sigma^{-b}(u(x,x))  \, v(x,x) \big\} \,.\,m\,.\, \big\{\sigma^{-b}(p(x,x)) \, \sigma^{-a-b}(q(x,x)) \big\}
	~ \in ~ {}_{a+b}M_{-a-b} \ .
$$
\epf

\subsection{A functor from $\TL_\kappa$ to $\bZ_d$-equivariant objects in $\Pd$}\label{sec:TL-to-PdZd}

We write $\Pd^{\bZ_d}$ for the category of $\bZ_d$-equivariant objects in $\Pd$, where the $\bZ_d$ action is given by the functor $A$ defined in the previous section. The definition and properties of categories of equivariant objects are recalled in Appendix\ref{ceo}.

By Theorem \ref{upn}, our aim now is to find a tensor functor
$$
	F : \TL_\kappa \to  \Pd^{\bZ_d} \ .
$$ 
According to Section \ref{uptl}, to construct a functor out of $\TL_\kappa$, we need to give a self dual object, duality maps, and compute the resulting constant $\kappa$. We will proceed as follows:
\begin{enumerate}
\item Give a self dual object $T \in \Pd$.
\item Give duality maps $u,n$, show they satisfy the zig-zag identities \eqref{eq:zigzag}, and compute $\kappa$.
\item Put a $\bZ_d$-equivariant structure on $T$ and show that the maps $u,n$ are $\bZ_d$-equivariant.
\end{enumerate}

\noindent
{\em Step 1:}
We listed self-dual objects of the from $P_{a:\lambda}$ at the end of Section \ref{sec:permutation-mf}. 
By Theorem \ref{thm:br-decomp}, there are only two choices which match the tensor products required by Corollary \ref{ftl}. In both cases, $d$ is odd, and either $\lambda=1$, $a = (d-1)/2$ or $\lambda=d-3$, $a = (d-3)(d-1)/2$. Both choices can be used in the construction below; we will work with the first option:
$$
	\text{$d$ odd}
	\quad , \qquad T := P_{\frac{d-1}2:1} = P_{\{\frac{d-1}2,\frac{d+1}2\}} \ .
$$	
Explicitly,
$$
T ~:\quad \xymatrix{\bC[x,y] \ar@/^10pt/[rrr]^{K(x,y)} &&& \bC[x,y] \ar@/^10pt/[lll]^{\frac{x^d-y^d}{K(x,y)}}}\quad ,
$$
where
$$
K(x,y)=\left( x-\eta^{\frac{d-1}{2}}y \right)\left( x-\eta^{\frac{d+1}{2}}y \right)=x^2+y^2+\kappa xy,\qquad \kappa = -(\eta^{\frac{d-1}{2}} + \eta^{\frac{d+1}{2}}) = 2 \cos  \tfrac{\pi}{d}  \ .
$$
Writing $\kappa$ for the coefficient of $xy$ will be justified below, where we will find it to be the parameter in $\TL_\kappa$.

\medskip
\noindent
{\em Step 2:}
Denote the isomorphism given in \eqref{eq:PS+P-S-iso} by $t : T \to T^+$, $t = (\id,-\id)$. Define maps $u : T \otimes T \to I$ and $n : I \to T \otimes T$ via
\be\label{eq:u-n-def}
	u = \big[ T \otimes T \xrightarrow{t \otimes \id} T^+ \otimes T \xrightarrow{ev_T} I  \big]
	\quad , \quad
	n = \big[ I \xrightarrow{coev_T} T \otimes T^+ \xrightarrow{\id \otimes t^{-1}} T \otimes T \big] \ .
\ee
{}From this one computes $u \circ n = \kappa$. For example,
$$
	u_0 \circ n_0 = A_T(x+z+\kappa y) = \kappa \ .
$$
Together with the zig-zag identities for $ev_T$ and $coev_T$ established in \cite{cr2} we have proved:

\bpr
$u$ and $n$ are morphisms in $\ZMFbi$. The satisfy the zig-zag identities in $\HMFbi$, as well as $n \circ u = \kappa$.
\epr

\noindent
{\em Step 3:}
We can make the $P_S$ $\bZ_d$-equivariant via
\be\label{eq:PS-Zd-equiv-structure}
	\tau_{S;a} : P_S \to {}_{a}(P_S)_{-a}
	\quad , \qquad
	\tau_{S;a} = \eta^{\frac{d+1}2 \,a (|S|-1)} \, s_{a,-a} \ ,
\ee
where $s_{a,-a}$ was given in \eqref{eq:PS-with-twisted-action-iso-to-PS}. 
These maps satisfy ${}_a(\tau_{S;b})_{-a} \circ \tau_{S;a} = \tau_{S;a+b}$, as required (cf.\ Appendix \ref{ceo}).
Note that on $I=P_{\{0\}}$, the above $\bZ_d$-equivariant structure is just $s_{a,-a} : I \to {}_aI_{-a}$, in agreement with the one on the tensor unit of $\Pd^{\bZ_d}$ as prescribed by Lemma \ref{lem:a()b-functor} and Proposition \ref{prop:CG-is-tensor}.

\ble\label{lem:ev-coev-t-Zd-equiv}
The maps $ev_{P_S}$ and $coev_{P_S}$ composed with the isomorphism $P_{-S} \simeq (P_S)^+$ from \eqref{eq:PS+P-S-iso} are $\bZ_d$-equivariant.
\ele

\bpf
For $coev$ we need to check commutativity of
$$
\xymatrix{ I  \ar[rr]^-{coev_{P_S}} \ar[d]_-{s_{a,-a}} && P_S \otimes (P_S)^+ \ar[r]^-{\sim} & P_S \otimes P_{-S} \ar[d]^-{\tau_{S;a} \otimes \tau_{-S;a}}
\\
{}_aI_{-a} \ar[rr]^-{{}_a(coev_{P_S})_{-a}} && {}_a(P_S \otimes (P_S)^+)_{-a} \ar[r]^-{\sim} & 
\begin{minipage}{9em} 
${}_a(P_S)_{-a} \otimes {}_a(P_{-S})_{-a}$ \newline
$=\, {}_a(P_S \otimes P_{-S})_{-a}$   
\end{minipage}
}
$$
which follows straightforwardly by composing the various maps. The corresponding diagram for $ev$ is checked analogously.
\epf

\bco
$u$ and $n$ are $\bZ_d$-equivariant morphisms.
\eco

According to Section \ref{sec:univ},
at this point we proved the existence of the tensor functor $\TL_\kappa \to  \Pd^{\bZ_d}$. To describe its image and to show that it annihilates the non-trivial tensor ideal in $\TL_\kappa$, we need to introduce a graded version of the above construction.

\subsection{Graded matrix factorisations}

There are several variants of graded matrix factorisations, see e.g.\ \cite{khovroz,Hori:2004ja,Wu09,cr1}.
The following one is convenient for our purpose. 
We take the grading group to be $\bC$, which  is natural from the relation to the R-charge in conformal field theory, but other groups are equally possible. For example, to construct the tensor equivalence in Theorem \ref{thm:main}, the grading group $d^{-1} \bZ$ is sufficient.

\bde
Let $S$ be a $\bC$-graded $\kk$-algebra such that $W \in S$ has degree 2. A {\em $\bC$-graded matrix factorisation} of $W$ over $S$ is a matrix factorisation $(M,d)$ of $W$ over $S$ such that the $S$ action on $M$ is compatible with the $\bC$-grading and $d$ has $\bC$-degree 1. That is, if $q(s)$ (resp.\ $q(m)$) denotes the $\bC$-degree of a homogeneous element of $S$ (resp.\ $M$), then $q(s.m) = q(s)+q(m)$ and $q(d(m)) = q(m)+1$.
\ede

In analogy with Section \ref{sec:cat-mf} we define 
$\MFgr_{S,W}$, $\ZMFgr_{S,W}$ and $\HMFgr_{S,W}$
to have $\bC$-graded matrix factorisations as objects and only $\bC$-degree zero morphisms. For example, 
\begin{align*}
	\HMFgr_{S,W}(M,N) ~=~ & \big\{ f \in \ZMF_{S,W}(M,N) \,| \text{ $f$ has $\bC$-degree $0$ }\big\} 
	\nonumber \\
	& ~ /~ \big\{  \delta(g)\, \big| \text{ $g : M \to N$ is $S$-linear, $\bZ_2$-odd and of $\bC$-degree $-1$ } \big\} \ .
	\nonumber
\end{align*}
The same definitions apply to matrix bifactorisations, giving categories $\MFgr_{\mathrm{bi};S,W}$, etc. Under tensor products, the $\bC$-degree is additive. 

\medskip

We will again restrict our attention to the case $S = \bC[x]$ and $W = x^d$, so that $q(x) = \frac2d$.

\medskip

As an example, let us describe all $\bC$-gradings on permutation type matrix bifactorisations. The $\bC$-grading on $\bC[x,y]$ is fixed by choosing the degree of 1. Let thus $\bC[x,y]\{\alpha\}$ be the graded $\bC[x]$-$\bC[x]$-bimodule with $q(1) = \alpha$. The possible $\bC$-gradings on $P_S$ are
$$
P_{S}\{\alpha\} ~:\quad\xymatrix{\bC[x,y]\{\alpha+\frac2d|S| -1\}  \ar@/^10pt/[rrr]^{d_1 = \prod\limits_{j \in S}(x-\eta^j y)} &&& \bC[x,y]\{\alpha\}\ar@/^10pt/[lll]^{d_0 = \prod\limits_{j \in \overline S}( x-\eta^j y)}} \quad .
$$
The unit isomorphism $\lambda_M$ given in \eqref{eq:unit-isos} above becomes a morphism in $\HMFbigr$ precisely if the unit object is $\bC$-graded as
$$
	I = P_{\{0\}}\{0\} \ .
$$
To see this note that $x^m y^n \in I_0 = \bC[x,y]$ will act as a degree $2(m+n)/d$ map on $M$. With this charge assignment for $I$, $\HMFbigr$ is tensor.

Next we work out the grading on $M^+$ for $M$ with $M_0$ and $M_1$ of rank 1. One first convinces oneself that for a homogeneous $p \in \bC[x,y,z]$ we have $\mathrm{deg}(A_M(p)) = \mathrm{deg}(p) - \mathrm{deg}(d_1^M(x,y))+1$, where $\mathrm{deg}$ denotes the polynomial degree. So if $M_0 = \bC[x,y]\{\alpha\}$, for $A_M$ to give a $\bC$-degree 0 map, we need $M^+_0 = \bC[x,y]\{-\alpha+\frac2d(1-\mathrm{deg}(d^M_1))\}$ (cf.\ \cite[Sect.\,2.2.4]{cr2}). This forces the $\bC$-grading to be
\begin{align*}
M~&:~
\xymatrix{\bC[x,y]\{\alpha+\frac2d\mathrm{deg}(d_1)-1\} \ar@/^10pt/[rr]^{d_1(x,y)} && \bC[x,y]\{\alpha\} 
\ar@/^10pt/[ll]^{d_0(x,y)}}
\nonumber\\
\leadsto\ M^+~&:~
\xymatrix{\bC[x,y]\{-\alpha-1+\frac2d\} \ar@/^10pt/[rr]^{d_1^{M^+} := -d_1(y,x)} &&\bC[x,y]\{-\alpha+\frac2d(1-\mathrm{deg}(d_1))\} \ar@/^10pt/[ll]^{d_0^{M^+} := d_0(y,x)}}\quad .
\end{align*}
One can check that $ev$ and $coev$ are indeed degree 0 maps with respect to these gradings. Note that we have $I^+=I$ also as graded matrix bifactorisations.

In the next section we will be interested in the $P_S\{\alpha\}$ with $\alpha = \frac{1-|S|}d$. We abbreviate these as $\hat P_S$. This subset of graded permutation type matrix bifactorisations 
is closed under taking duals:
$$
	(\hat P_S)^+ \simeq \hat P_{-S}
	\qquad , \quad \text{where} \quad \hat P_S = P_S\big\{\tfrac{1-|S|}d\big\} \ .
$$
An explicit isomorphism is again given by \eqref{eq:PS+P-S-iso}, which is easily checked to have $\bC$-degree 0.

The next two lemmas show that the $\hat P_S$ generate (under direct sums) a semi-simple subcategory of $\HMFbigr$.

\ble
$\ZMFbigr(\hat P_{R}, \hat P_{S})$ is $\bC \id$ if $R=S$ and $0$ else.
\ele

\bpf
Write $\alpha = \frac{1-|R|}d$ and $\beta = \frac{1-|S|}d$, such that $\hat P_R = P_R\{\alpha\}$ and $\hat P_S = P_S\{\beta\}$. 
The morphism space $\ZMFbi(P_R,P_S)$ is given by all $(p,q)$ with $p,q \in \bC[x,y]$ such that $p \cdot d^{P_R}_1 = d^{P_S}_1 \cdot q$. For $(p,q)$ to be also in $\ZMFbigr(P_R\{\alpha\},P_S\{\beta\})$, we need  $p,q$ to be homogeneous and $\alpha = \beta + \frac2d\deg(p)$ and $\alpha + \frac2d |R|-1 = \beta + \frac2d |S|-1 + \frac2d\deg(q)$. This simplifies to $2 \deg(p) = |S|-|R|$ and $2 \deg(q) = |R|-|S|$, which is possible only for $|R|=|S|$, in which case $p,q$ are constants.
Finally, the condition $p \cdot d^{P_R}_1 = d^{P_S}_1 \cdot q$ has non-zero constant solutions only if $R=S$.
\epf

\ble \label{lem:hmfgrPSPR}
$\hat P_\emptyset$ and $\hat P_{\bZ_d}$ are zero objects in $\HMFbigr$.
For $R,S \neq \emptyset,\bZ_d$ we have $\HMFbigr(\hat P_{R}, \hat P_{S})=\ZMFbigr(\hat P_{R}, \hat P_{S})$.
\ele

\bpf
That $\hat P_\emptyset$ and $\hat P_{\bZ_d}$ are zero objects in $\HMFbigr$ follows since one component of the twisted differential is $1$, and hence there is a contracting homotopy for the identity morphism. 

Let now $R,S$ be nonempty proper subsets of $\bZ_d$.
For the second part of the statement one checks that there are no $\bZ_2$-odd morphisms of $\bC$-degree $-1$ from $\hat P_R$ to $\hat P_S$. For example, a $\bC$-degree $-1$ map $\psi_0 : (\hat P_{R})_0 \to (\hat P_{S})_1$ has to satisfy
$$
  \frac{1+|S|}{d}-1+\frac{2\,\mathrm{deg}(\psi_0(x,y))}{d} - \frac{1-|R|}{d} = -1 \ ,
$$
where $\mathrm{deg}(\psi_0)$ is the polynomial degree of $\psi_0(x,y)$. Thus $\mathrm{deg}(\psi_0) =- \frac{|S|+|R|}2$, and $\psi_0$ can be non-zero only if $|R|=|S|=0$. An analogous computation for $\psi_1$ shows $\mathrm{deg}(\psi_1) = \frac{|S|+|R|}2 - d$, and so $\psi_1$ can be non-zero only if $|R|=|S|=d$.
\epf

We now focus on the graded matrix factorisations $\hat P_{a:\lambda}$, i.e.\ the $P_S\{\alpha\}$ with $S = \{a,a{+}1,\dots,a{+}\lambda\}$ and $\alpha = - \lambda/d$. We define
$$
	\Pdgr = \big\langle \hat P_{a:\lambda} \,\big|\, a \in \bZ_d, \lambda \in \{0,\dots,d-2\} \big
	\rangle_\oplus
	~\subset~ \HMFbigr \ ,
$$
i.e.\ the full subcategory of $\HMFbigr$ consisting of objects isomorphic, in $\HMFbigr$, to finite direct sums of the $\hat P_{a:\lambda}$.

We now need to check whether the decomposition of tensor products in Theorem \ref{thm:br-decomp} carries over to the graded case. This could be done by adapting the method used in \cite{brunnerrogg}, which works in the stable category of $\bC[x,y]/\langle x^d-y^d \rangle$ modules. We give a related but different proof by providing explicit $\bC$-charge 0 embeddings of the direct summands in the decomposition of $\hat P_{a:1} \otimes \hat P_{b:\lambda}$ and proving that they give an isomorphism via Remark \ref{rem:H(M)}. This is done in Appendix \ref{app:graded-tensor}

\bth\label{thm:Pdgr-tensor-closed}
The category $\Pdgr$ is semi-simple with simple objects $\hat P_{a:\lambda}$, $a \in \bZ_d$ and $\lambda \in \{0,1,\dots,d-2\}$. It is closed under tensor products and the direct sum decomposition of $\hat P_{m:\lambda} \otimes \hat P_{n:\nu}$ in $\HMFbigr$ is as in Theorem \ref{thm:br-decomp}.
\eth

\subsection{A functor from $\TL_\kappa$ to $\bZ_d$-equivariant objects in $\Pdgr$}

The morphisms $\mu_{a,b}$ in \eqref{eq:mu-morph} have $\bC$-degree 0. The functor $\chi$ in Proposition \ref{prop:chi-functor} therefore also defines a tensor functor
$$
	\chi : \underline{\bZ_d} \longrightarrow \Pdgr \ .
$$
As in Section \ref{sec:Zd_to_Pd} we obtain two tensor functors $A , \Ad_\chi : \underline{\bZ_d} \to \mathrm{Aut}_\otimes(\Pdgr)$. The natural monoidal isomorphism $A \to \Ad_\chi$ established in Lemma \ref{lem:A-iso-Ad} uses only $\bC$-degree 0 morphisms. 

\medskip

Next we follow the three steps in Section \ref{sec:TL-to-PdZd} and verify that they carry over to the $\bC$-graded setting. Consider the self-dual object $\hat T \in \Pdgr$. The duality maps $n$ and $u$ from \eqref{eq:u-n-def} are of $\bC$-degree 0 since $t$, $ev_T$, $coev_T$ are. The maps $\tau$ from \eqref{eq:PS-Zd-equiv-structure} are equally of $\bC$-degree 0 and hence equip $\hat T$ with a $\bZ_d$-equivariant structure. The proof of Lemma \ref{lem:ev-coev-t-Zd-equiv} still applies and shows that $u$ and $n$ are $\bZ_d$-equivariant morphisms in $\HMFbigr$.

\medskip

By Section \ref{uptl} the data $\hat T$, $\tau$, $u$ and $n$ determine a tensor functor
$$
	F : \TL_\kappa \to  (\Pdgr)^{\bZ_d} \ .
$$ 
Here we used that $\hat T \in \Pdgr$ and that by Theorem \ref{thm:Pdgr-tensor-closed}, $\Pdgr$ is a full tensor subcategory of $\HMFbigr$.

\bth\label{thm:main}
There is a tensor equivalence $G : \C(N{=}2,d)_{NS} \to \Pdgr$ such that $G([l,l+2m]) \simeq \hat P_{m:l}$.
\eth

\bpf 
Corollary \ref{ftl} and the tensor product established in Theorem \ref{thm:Pdgr-tensor-closed} show that $F$ is not faithful and induces a fully faithful embedding $\tilde F:\T_\kappa \to  (\Pdgr)^{\bZ_d}$. By Theorem \ref{upn} the embedding $\tilde F$ gives rise to the functor $G: \C(N{=}2,d)_{NS}\to \Pdgr$.
The functor $G$ is fully faithful (it sends simple objects to simple objects) and surjective on (simple) objects.
Thus $G$ is an equivalence.

Recall that the $\bZ_d$-action on $\Pdgr$ is such that $a \in \bZ_d$ gets mapped to ${}_aI \cong P_{\{-a\}}$, and that $\tilde F$ maps $T \in \T_\kappa$ to $\hat P_{\frac{d-1}2:1} \in \Pdgr$. We choose the monoidal embedding $\underline \bZ_d \to \C(N{=}2,d)$ as $a \mapsto [0,-2a]$ (to avoid this minus sign, one can define $\chi$ in Proposition \ref{prop:chi-functor} as $\chi(a) = {}_{-a}I$, resulting in lots of minus signs in other places). The induced tensor functor $G$ obeys $G([1,d]) = \hat P_{\frac{d-1}2:1}$ and $G([0,2a]) = P_{\{a\}}$.
\epf

\bre\label{rem:galois}
Note that one can replace $\eta$ with any other primitive $d$'th root
of unity $\eta^l$ (here $l$ is coprime to $d$).
In particular replacing $\eta$ with $\eta^l$ in \eqref{pfs} gives another
matrix bifactorisation,
$P_S(\eta^l)$.
It is not hard to see that $P_{\{\frac{d-1}2,\frac{d+1}2\}}(\eta^l)$ is a
self-dual object of dimension $\kappa_l=2 \cos\tfrac{\pi l}{d}$
and defines a fully faithful embedding $\T_{\kappa_l}\to \HMFbi$.
Its image is additively generated by the direct summands in tensor powers of $P_{\{\frac{d-1}2,\frac{d+1}2\}}(\eta^l)$ and can be computed explicitly from Theorem \ref{thm:Pdgr-tensor-closed} with $\eta^l$ in place of $\eta$. This is an instance of the action of a Galois group on categories of matrix factorisation, see \cite[Rem.\,2.9]{crcr} for a related discussion.
\ere

\begin{appendices}

\section{Equivariant objects and
pointed categories}\label{sec:equiv+pointed}

Here we collect some (well known) categorical trivialities which allow us to avoid difficult calculations with matrix bifactorisations. 
Throughout Appendix \ref{sec:equiv+pointed} all tensor (and in particular all fusion) categories will be assumed to be strict. 
In labels for some arrows in our diagrams we suppress tensor product symbols for compactness. 

\subsection{Categories of equivariant objects}\lb{ceo}

Let $G$ be a group. 
An {\em action} of $G$ on a tensor category $\C$ is a monoidal functor $F:\underline{G}\to\Aut_\ot(\C)$ from the discrete monoidal category $\underline{G}$ to the groupoid of tensor autoequivalences of $\C$.
More explicitly, a $G$-action on $\C$ consists of a collection $\{F_g\}_{g\in G}$ of tensor autoequivalences $F_g:\C\to\C$ labelled by elements of $G$ together with natural isomorphisms $\phi_{f,g}:F_f\circ F_g\to F_{fg}$ of tensor functors such that $\phi_{f,e} = 1,\phi_{e,g} =1$ and such that the diagram 
$$\xymatrix{F_f\circ F_g\circ F_h\ar[rr]^{\phi_{f,g}\circ 1} \ar[d]_{1\circ\phi_{g,h}} && F_{fg}\circ F_g \ar[d]^{\phi_{fg,h}} \\ 
F_f\circ F_{gh} \ar[rr]^{\phi_{f,gh}} && F_{fgh}}$$
commutes for any $f,g,h\in G$.

Let $\C$ be a tensor category together with a $G$-action.
An object $X\in\C$ is {\em $G$-equivariant} if it comes equipped with a collection of isomorphisms $x_g:X\to F_g(X)$ such that 
the diagram 
$$\xymatrix{X\ar[rr]^{x_{fg}} \ar[d]_{x_{f}} && F_{fg}(X) \\ 
F_f(X) \ar[rr]_{F_f(x_g)} && F_f(F_g(X)) \ar[u]_{\phi_{f,g}} }$$
commutes for any $f,g\in G$.

A morphism $a:X\to Y$ between $G$-equivariant objects $(X,x_g), (Y,y_g)$ is {\em $G$-equivariant} if 
the diagram 
$$\xymatrix{X\ar[rr]^{x_g} \ar[d]_{a} && F_{g}(X) \ar[d]^{F_g(a)}\\ 
Y \ar[rr]^{y_g} && F_g(Y) }$$
commutes for any $g\in G$.
Denote by $\C^G$ the category of $G$-equivariant objects in $\C$.

\bpr\label{prop:CG-is-tensor}
Let $\C$ be a strict tensor category with a $G$-action. 
Then the category $\C^G$ is strict tensor with 
tensor product $(X,x_g)\ot(Y,y_g) = (X\ot Y,(x|y)_g)$, where $(x|y)_g$ is defined by
$$\xymatrix{X\ot Y \ar[rr]^(.4){x_gy_g} && F_g(X)\ot F_g(Y) \ar[rr]^(.55){(F_g)_{X,Y}} && F_g(X\ot Y)}$$
and with unit object $(I,\iota)$, where $\iota_g:I\to F_g(I)$ is the unit isomorphism of the tensor functor $F_g$.
\epr
\bpf
All we need to check is that the $G$-equivariant structures of the triple tensor products $(X,x_g)\ot((Y,y_g)\ot(Z,z_g))$ and $((X,x_g)\ot(Y,y_g))\ot(Z,z_g)$ coincide. These $G$-equivariant structures $x|(y|z), (x|y)|z$ are the top and the bottom outer paths of the diagram
$$\xymatrix{X\ot Y\ot Z\ar[rrd]^{x_gy_gz_g} \\ && F_g(X)\ot F_g(Y)\ot F_g(Z) \ar[rr]^{(F_g)_{X,Y}\id} \ar[d]_{\id(F_g)_{Y,Z}} && F_g(X\ot Y)\ot F_g(Y) \ar[d]^{(F_g)_{XY,Z}} \\ &&
F_g(X)\ot F_g(Y\ot Z) \ar[rr]^{(F_g)_{X,YZ}} && F_g(X\ot Y\ot Z)}$$
whose commutativity is the coherence of the tensor structure of $F_g$.
\epf

Clearly the forgetful functor
$$\C^G\ \to\ \C,\qquad (X,x)\mapsto X$$ is tensor.

\bre\lb{rem2}
It is possible to define more general $G$-actions on (tensor) categories involving associators for $G$ (3-cocycles for $G$).
All constructions generalise straightforwardly.
\ere

\subsection{Inner actions and monoidal centralisers of pointed subcategories}\lb{iac}

An object $P$ of a tensor category $\C$ is {\em invertible} if the dual object $P^*$ exists and the evaluation $ev_P:P^*\ot P\to I$ and coevaluation $coev_P:I\to P\ot P^*$ maps are isomorphisms. Clearly an invertible object is simple since $\C(P,P)\simeq\C(I,I)=\kk$.

The set $Pic(\C)$ of isomorphism classes of invertible objects is a group with respect to the tensor product (the {\em Picard group} of $\C$).
Choosing a representative $s(p)$ in each isomorphism class $p\in Pic(\C)$ and isomorphisms $\sigma(p,q):s(p)\ot s(q)\to s(pq)$ for each pair $p,q\in Pic(\C)$ allows us to define a function $\alpha:Pic(\C)^{\times 3}\to \kk^*$ (here $\kk^*$ is the multiplicative group of non-zero elements of $\kk$). Indeed for $p,q,r\in Pic(\C)$ the composition 
$$
s(pqr) \xrightarrow{\sigma(pq,r)^{-1}}  s(pq)\ot s(r) \xrightarrow{\sigma(p,q)^{-1}\,\id}  s(p)\ot s(q)\ot s(r) \xrightarrow{\id\,\sigma(q,r)}  s(p)\ot s(qr) \xrightarrow{\sigma(p,qr)}  s(pqr)$$
is an automorphism of $s(pqr)$ and thus has a form $\alpha(p,q,r)\id_{s(pqr)}$ for some $\alpha(p,q,r)\in \kk^*$.
It is easy to see that $\alpha$ is a 3-cocycle and that the class $[\alpha]\in H^3(Pic(\C),\kk^*)$ does not depend on the choice of $s$ and $\sigma$. 

A tensor category $\C$ is {\em pointed} if all its simple objects are invertible. 
A fusion pointed category $\C$ can be identified with the category $\V(G,\alpha)$ of $G$-graded vector spaces, where $G=Pic(\C)$ and with the associativity constraint twisted by $\alpha\in H^3(Pic(\C),\kk^*)$.

Let $P$ be an invertible object of a tensor category $\C$.
The functor
$$P\ot-\ot P^*:\C\to \C,\qquad X\mapsto P\ot X\ot P^*$$ comes equipped with a monoidal structure 
$$\xymatrix{P\ot X\ot P^*\ot P\ot Y\ot P^* \ar[rr]^(.6){\id\,ev_P\,\id} && P\ot X\ot Y\ot P^*}$$
making it a tensor autoequivalence, the {\em inner autoequivalence} corresponding to $P$.
The assignment $P\mapsto P\ot-\ot P^*$ defines a homomorphism of groups
$Pic(\C)\ \to\ Aut_\ot(\C)\ .$

The {\em monoidal centraliser} $\Z_\D(F)$ of a tensor functor $\C\to\D$ is the category of pairs $(Z,z)$, where $Z\in\D$ and $z_X:Z\ot F(X)\to F(X)\ot Z$ are a collection of isomorphisms, natural in $X\in\C$, such that $Z_I=\id$ and such that the diagram
$$\xymatrix{Z\ot F(X\ot Y) \ar[d]_{\id\,F_{X,Y}} \ar[rr]^{z_{XY}} && F(X\ot Y)\ot Z \ar[d]^{F_{X,Y}\,\id} \\ Z\ot F(X)\ot F(Y)\ar[dr]_{z_X\,\id} && F(X)\ot F(Y)\ot Z\\ & F(X)\ot Z\ot F(Y) \ar[ru]_{\id\,z_Y} }$$
commutes for any $X,Y\in\C$.
A morphism $(Z,z)\to(Z',z')$ in $\Z_\D(F)$ 
is a morphism $f:Z\to Z'$ in $\D$ such that the diagram
$$\xymatrix{Z\ot F(X)\ar[rr]^{z_X} \ar[d]_{f\,\id} && F(X)\ot Z \ar[d]^{\id\,f} \\ Z'\ot F(X)\ar[rr]^{z'_X} && F(X)\ot Z'}$$
commutes for any $X\in\C$.

\bpr
Let $F:\C\to\D$ be a tensor functor between strict tensor categories. 
Then the monoidal centraliser $\Z_\D(F)$ is strict tensor with the tensor product $(Z,z)\ot(Z',z') = (Z\ot Z',z|z')$ where $(z|z')_X$ is defined by
$$\xymatrix{Z\ot Z'\ot F(X)\ar[rr]^{(z|z')_{X}} \ar[dr]_{\id\,z'_X} &&
F(X) \ot Z\ot Z' 
 \\ & Z\ot F(X)\ot Z' \ar[ru]_{z_X\,\id} }$$
and with the unit object $(I,1)$.
\epr
\bpf
Note that the monoidal centraliser $\Z_\D(Id_\D)$ of the identity functor $Id_\D:\D\to\D$ is the monoidal centre $\Z(\D)$. 
The proof of the proposition is identical to the proof of monoidality of the monoidal centre (see \cite{js0}).
\epf

Clearly the forgetful functor
$$\Z_\D(F)\ \to\ \D,\qquad (Z,z)\mapsto Z$$ is tensor.

Let $G$ be a group and $\V(G)$ be the pointed tensor category whose group of isomorphism classes of objects is $G$ and which has trivial associator.
A tensor functor $F:\V(G)\to\C$ gives rise to the action of $G$ on $\C$ by inner autoequivalences $F_g(X) = F(g)\ot X\ot F(g)^*,\ g\in G$.

\bth\lb{cps}
Let $G$ be a group and let $\C$ be a tensor category with a tensor functor $F:\V(G)\to\C$. 
Then the monoidal centraliser $\Z_\C(F)$ is tensor equivalent to the category of $G$-equivariant objects $\C^G$, where the $G$-action is defined by the functor $\V(G)\to\C$ as above.
\eth
\bpf
Define a functor $\Z_\C(F)\to \C^G$ by assigning to $(Z,z)\in\Z_\C(F)$ a $G$-equivariant object $(Z,\tilde z_g)_{g\in G}$ with $\tilde z_g:Z\to F_g(X) = F(g)\ot X\ot F(g)^*$ given by
$$\xymatrix{Z \ar[rr]^-{\id\,coev_{F(g)}} && Z\ot F(g)\ot F(g)^*\ar[rr]^{z_g\,\id} && F(g)\ot Z\ot F(g)^*}\quad .$$
It is straightforward to see that this is a tensor equivalence.
\epf

\subsection{Tensor functors from products with pointed categories}

Recall from \cite{de,baki}
that the {\em Deligne product} $\C\boxtimes\D$ of $\kk$-linear semi-simple categories $\C$ and $\D$ is a semi-simple category with simple objects $X\boxtimes Y$ for $X$ and $Y$ being simple objects of $\C$ and $\D$ correspondingly.
One can extend the definition of $X\boxtimes Y$ to arbitrary $X\in\C$ and $Y\in\D$. The hom spaces between these objects are
$$(\C\boxtimes\D)(X\boxtimes Y,X'\boxtimes Y') = \C(X,X')\ot_\kk \D(Y,Y')\ ,$$
where on the right is the tensor product of vector spaces over $\kk$.

The Deligne product of fusion categories is fusion with the unit object $I\boxtimes I$ and the tensor product defined by
$$(X\boxtimes Y)\ot(X'\boxtimes Y') = (X\ot X')\boxtimes (Y\ot Y')\ .$$
The Deligne product of fusion categories has another universal property, which we describe next.

We say that a pair of tensor functors $F_i:\C_i\to\D$ has {\em commuting images} if they come equipped with a collection of isomorphisms
$c_{X_1,X_2}:F_1(X_1)\ot F_2(X_2)\to F_2(X_2)\ot F_1(X_1)$ natural in $X_i\in\C_i$ and such that the following diagrams commute for all $X_i,Y_i\in \C_i$:
$${\xymatrixcolsep{1pc}\xymatrix{F_1(X_1)\ot F_2(I)\ar[rr]^{c_{X_1,I}} \ar[d] && F_2(I)\ot F_1(X_1) \ar[d] \\ F_1(X_1)\ot I\ar@{=}[r] & F_1(X_1) \ar@{=}[r] & I\ot F_1(X_1)}}\qquad 
{\xymatrixcolsep{1pc}\xymatrix{F_1(I)\ot F_2(X_2)\ar[rr]^{c_{I,X_2}} \ar[d] && F_2(X_2)\ot F_1(I) \ar[d] \\ I\ot F_2(X_2)\ar@{=}[r] & F_2(X_2) \ar@{=}[r] & F_2(X_2)\ot I}}$$
$$\xymatrix{F_1(X_1\ot Y_1)\ot F_2(X_2) \ar[rr]^{c_{X_1Y_1,X_2}} \ar[d]_{(F_1)_{X_1,Y_1}1} && F_2(X_2)\ot F_1(X_1\ot Y_1) \ar[d]^{1(F_1)_{X_1,Y_1}} \\ F_1(X_1)\ot F_1(Y_1)\ot F_2(X_2)  \ar[rd]_{1c_{Y_1,X_2}} && F_2(X_2)\ot F_1(X_1)\ot F_1(Y_1) \\ & F_1(X_1)\ot F_2(X_2)\ot F_1(Y_1)  \ar[ru]_{c_{X_1,X_2}1}}$$
$$\xymatrix{F_1(X_1)\ot F_2(X_2\ot Y_2) \ar[rr]^{c_{X_1,X_2Y_2}} \ar[d]_{1(F_2)_{X_2,Y_2}} && F_2(X_2\ot Y_2)\ot F_1(X_1) \ar[d]^{(F_2)_{X_2,Y_2}1} \\ F_1(X_1)\ot F_2(X_2)\ot F_2(Y_2)  \ar[rd]_{c_{X_1,X_2}1} && F_2(X_2)\ot F_2(Y_2)\ot F_1(X_1) \\ & F_2(X_2)\ot F_1(X_1)\ot F_2(Y_2)  \ar[ru]_{1c_{X_1,Y_2}}}$$

\bpr\lb{upd}
The Deligne product $\C_1\boxtimes\C_2$ of fusion categories $\C_1$ and $\C_2$ is the initial object among pairs of tensor functors $F_i:\C_i\to\D$ with commuting images, that is for a pair of tensor functors $F_i:\C_i\to\D$ with commuting images there is a 
	unique
tensor functor $F:\C_1\boxtimes\C_2\to\D$ making the diagram 
$$\xymatrix{\C_1 \ar[rd] \ar[rdd]_{F_1} && \C_2 \ar[ld] \ar[ldd]^{F_2} \\ & \C_1\boxtimes\C_2 \ar@{.>}[d]^F \\ & \D}$$
commutative.
\epr
\bpf
Note that the assignments $X_1\mapsto X_1\boxtimes I,\ X_2\mapsto I\boxtimes X_2$ define a pair of tensor functors $\C_i\to C_1\boxtimes\C_2$ with commuting images.

Conversely let $F_i:\C_i\to\D$ be a pair of tensor functors with commuting images.
Define $F:\C_1\boxtimes\C_2\to\D$ by $F(X_1\boxtimes X_2) = F_1(X_1)\ot F_2(X_2)$. 
	Since $\C_1$ and $\C_2$ are fusion, this determines $F$ uniquely as a $\kk$-linear functor.
	The monoidal structure for $F$ is uniquely determined to be 
$$\xymatrix{
F(X_1\boxtimes X_2)\ot F(Y_1\boxtimes Y_2) \ar@{=}[d] \ar[rrr]^{F_{X_1\boxtimes X_2,Y_1\boxtimes Y_2}} &&& F\big((X_1\boxtimes X_2)\ot(Y_1\boxtimes Y_2)\big) \ar@{=}[d] \\
F(X_1)\ot F(X_2)\ot F(Y_1)\ot F(Y_2) \ar[d]_{1c_{X_2,Y_1}1} &&& F\big((X_1\ot Y_1)\boxtimes(X_2\ot Y_2)\big) \ar@{=}[d] \\
F(X_1)\ot F(Y_1)\ot F(X_2)\ot F(Y_2) \ar[rrr]^{(F_1)_{X_1,Y_1}(F_2)_{X_2,Y_2}} &&& F_1(X_1\ot Y_1)\ot F_2(X_2\ot Y_2) }$$
It is straightforward to check that this definition satisfies the coherence axioms of a monoidal structure. 
\epf

\bre\lb{ciz}
Note that the data of a pair of tensor functors $F_i:\C_i\to\D$ with commuting images amounts to a tensor functor $\C_1\to \Z_\D(F_2)$ whose composition with the forgetful functor $\Z_\D(F_2)\to\D$ equals $F_1$.
\ere

\bth\lb{fpp}
Let $\C$ be a fusion
category and let $G$ be a finite
group.
Then the data of a tensor functor $\C\boxtimes\V(G)\to\D$ amounts to a tensor functor $\V(G)\to\D$ and a tensor functor $\C\to \D^G$, where the $G$-action is defined by the functor $\V(G)\to\D$ as in Appendix \ref{iac}.
\eth
\bpf
By Proposition \ref{upd} a tensor functor $\C\boxtimes\V(G)\to\D$ corresponds to a pair of tensor functors $F:\V(G)\to\D,\ F':\C\to\D$ with commuting images.
By Remark \ref{ciz} this is equivalent to a tensor functor $\C\to\Z_\D(F)$ to the centraliser of $F$.
Finally by Theorem \ref{cps} the centraliser $\Z_\D(F)$ is canonically equivalent to the category of equivariant objects $\D^G$.
\epf

\bre
It is possible to extend Theorem \ref{fpp} to the case of pointed categories $\V(G,\alpha)$ with non-trivial associators $\alpha\in Z^3(G,\kk^*)$.
As for Remark \ref{rem2} all constructions generalise straightforwardly.
\ere

\section{Proof of Theorem \ref{thm:Pdgr-tensor-closed}}\label{app:graded-tensor}

Semi-simplicity of $\Pdgr$ follows from Lemma \ref{lem:hmfgrPSPR}, as does the list of simple objects.

\medskip

Let $\lambda,\mu \in \{0,1,\dots,d{-}2\}$ and $a,b \in \bZ_d$.
To show the decomposition rule 
\be\label{eq:P-hat-decomp}
	\hat P_{a:\lambda}  \otimes \hat P_{b:\mu} \simeq \bigoplus_{\nu=|\lambda-\mu|~\mathrm{step}\,2}^{\min(\lambda+\nu,2d-4-\lambda-\nu)} \hat P_{a+b-\frac12(\lambda+\mu-\nu):\nu} \quad ,
\ee
we verify the cases $\lambda=0$ and $\lambda=1$ explicitly. The general case follows by a standard argument using induction on $\lambda$.

\medskip
\noindent
{\em Case $\lambda=0$:} The isomorphism $\hat P_{a:0}  \otimes \hat P_{b:\mu} \simeq \hat P_{a+b:\mu}$ is immediate from the isomorphism $\hat P_{a:0} \simeq {}_{-a}I$ given in \eqref{eq:PS-with-twisted-action-iso-to-PS}, the isomorphism ${}_{-a}I \otimes \hat M \to {}_{-a}M$ provided by ${}_{-a}(\lambda_M)$ for any matrix factorisation $M$, and ${}_{-a}(P_{S}) \simeq P_{S+a}$, again from \eqref{eq:PS-with-twisted-action-iso-to-PS}.

\medskip
\noindent
{\em Case $\lambda=1$:}
For $\mu=0$ the isomorphism $\hat P_{a:1} \otimes \hat P_{b:0} \simeq P_{a+b:1}$ constructed as in case $\lambda=0$, using $P_{b:0} \simeq I_{-b}$. To show the decomposition \eqref{eq:P-hat-decomp} for $\mu \in \{1,2,\dots,d-2 \}$
we start by giving maps
$$
	g^- \,:\, \hat P_{a+b+1:\mu-1} \longrightarrow \hat P_{a:1} \otimes \hat P_{b:\mu}
	\quad , \quad
	g^+ \,:\, \hat P_{a+b:\mu+1} \longrightarrow \hat P_{a:1} \otimes \hat P_{b:\mu}
$$
in $\ZMFbigr$. Write $A = \hat P_{a:1}$, $B = P_{b:\mu}$, 
$Q_- = \hat P_{a+b+1:\mu-1}$ and $Q_+ = P_{a+b:\mu+1}$. We have to find $g^\varepsilon_{ij}$ that fit into the diagram
\be\label{eq:P1Plam-embed-summands}
\xymatrix{
Q_\ve \ar[ddd] & & \bC[x,z]\left\{\frac{\mu+2+\ve}{d}-1\right\} \ar@/^10pt/[rrrr]^{q_\ve(x,z)} \ar[ddd]_{\left(\begin{array}{c} \scriptstyle{g_{10}^\ve} \\ \scriptstyle{g_{01}^\ve} \end{array}\right)} &&&& \bC[x,z]\left\{-\frac{\mu+2}{d}\right\} \ar@/^10pt/[llll]^{\frac{x^d-z^d}{q_\ve(x,z)}} \ar[ddd]^{\left(\begin{array}{r} \scriptstyle{g_{00}^\ve} \\ \scriptstyle{g_{11}^\ve} \end{array}\right)} \\ \\ \\
A\ot B && {\begin{array}{c}{\bC[x,y,z]\left\{\frac{3-\mu}{d}-1\right\}} \\ {\op} \\ {\bC[x,y,z]\left\{\frac{1+\mu}{d}-1\right\}} \end{array} }\ar@/^10pt/[rrrr]^{\left(\begin{array}{rr} \scriptstyle{p_1(x,y)} & \scriptstyle{p_\mu(y,z)} \\ \frac{y^d-z^d}{p_\mu(y,z)} &  \frac{x^d-y^d}{p_1(x,y)} \end{array}\right)} &&&& 
{\begin{array}{c}{\bC[x,y,z]\left\{-\frac{\mu+1}{d}\right\}} \\ {\op} \\ {\bC[x,y,z]\left\{\frac{5+\mu}{d}-2\right\}} \end{array} }
\ar@/^10pt/[llll]^{\left(\begin{array}{rr} \frac{x^d-y^d}{p_1(x,y)} & \scriptstyle{-p_\mu(z,y)} \\ \frac{y^d-z^d}{p_\mu(y,z)} & \scriptstyle{p_1(x,y)} \end{array}\right)}
}
\ee
Here,
\begin{align*}
	p_1(x,y) &= (x-\eta^{a}y)(x-\eta^{a+1}y) ~,&
	q_-(x,z) &= \prod_{j=a+b+1}^{a+b+\mu} (x-\eta^j z) ~,
	\nonumber \\
	p_\mu(y,z) &= \prod_{j=b}^{b+\mu} (y-\eta^j z) ~, &
	q_+(x,z) &= \prod_{j=a+b}^{a+b+\mu+1} (x-\eta^j z) \ .
\end{align*}
Comparing $\bC$-degrees determines the polynomial degrees of the individual maps to be
\begin{align*}
\mathrm{deg}(g^\varepsilon_{10}) &= \mu-\tfrac12(1-\varepsilon)  ~,& 
\mathrm{deg}(g^\varepsilon_{00}) &= \tfrac12(1-\varepsilon) ~, 
\nonumber\\
\mathrm{deg}(g^\varepsilon_{01}) &= \tfrac12(1+\varepsilon) ~,&
\mathrm{deg}(g^\varepsilon_{11}) &= d-2 - \mu  -\tfrac12(1+\varepsilon) \ . 
\end{align*}
Commutativity of \eqref{eq:P1Plam-embed-summands} is equivalent to
\begin{align*}
&\text{(i)}& q_\varepsilon(x,z) \, g^\varepsilon_{00}(x,y,z) &= p_1(x,y) \, g^\varepsilon_{10}(x,y,z) + p_\mu(y,z) \, g^\varepsilon_{01}(x,y,z)
\nonumber \\
&\text{(ii)}& q_\varepsilon(x,z) \, g^\varepsilon_{11}(x,y,z) &= -\frac{y^d-z^d}{p_\mu(y,z)} \, g^\varepsilon_{10}(x,y,z) + \frac{x^d-y^d}{p_1(x,y)} \, g^\varepsilon_{01}(x,y,z)
\end{align*}
These conditions imply the remaining two conditions.
Let us show how one arrives at $g^-$ in some detail and then just state the result for $g^+$. 

\medskip

We have $\mathrm{deg}(g^-_{01})=0$ and we make the ansatz $g^-_{01}=1$ (choosing $g^-_{01}=0$ forces $g^-=0$, so this is really a normalisation condition).
The polynomial $g^-_{00}$ is of degree 1, so $g^-_{00}(x,y,z) = \alpha x + \beta y + \gamma z$ for some $\alpha,\beta,\gamma \in \bC$. Condition (i) determines $g^-_{10}$ uniquely to be
$$
	g^-_{10}(x,y,z) = \frac{q_-(x,z)(\alpha x + \beta y + \gamma z) - p_\mu(y,z)}{p_1(x,y)} \ .
$$
We need to impose the condition that $g^-_{10}$ is a polynomial. 
This amounts to verifying that the numerator has zeros for $y=\eta^{-a}x$ and $y=\eta^{-a-1}x$.
Using 
\begin{align*}
p_\mu(\mu^{-a}x,z) &=  q_-(x,z) \, \eta^{-a(\mu+1)} (x-\eta^{a+b}z) \ ,
\nonumber \\
p_\mu(\mu^{-a-1}x,z) &= q_-(x,z) \, \eta^{-(a+1)(\mu+1)}  (x-\eta^{a+b+\mu+1}z)
\end{align*}
gives the unique solution
$$
	g^-_{00}(x,y,z) = \eta^{-a\mu} \Bigg( - \eta^{-a-1} \frac{1-\eta^{-\mu}}{1-\eta^{-1}} \, x
	+ \frac{1 -\eta^{-\mu-1}}{1-\eta^{-1}} \, y
	- \eta^{b} \,z \Bigg) \ .
$$
Finally, a short calculation shows that condition (ii) is equivalent to
$$
	g_{11}^- = \frac{1}{p_1(x,y)}\Bigg( \frac{x^d-z^d}{q_-(x,z)} - \frac{y^d-z^d}{p_\mu(y,z)} \, g^-_{00}(x,y,z) \Bigg) \ .
$$
The term in brackets is clearly a polynomial. To show that it is divisible by $p_1(x,y)$, one simply verifies that the term is brackets is zero for $y = \eta^{-a}x$ and $y = \eta^{-a-1}x$.

For $g^+$ the calculation works along the same lines with the result
\begin{align*}
	g^+_{00} &= 1  
	\nonumber \\
	g^+_{01} &= 
		\eta^{a (\mu+1)  } \Bigg( \frac{1-\eta^{\mu+2}}{1-\eta} \, x 
		- \eta^{a+1} \frac{1 - \eta^{\mu+1}}{1-\eta} \, y - \eta^{a+b+\mu+1} \, z \Bigg)
	\nonumber \\
	g^+_{10} &= \frac{q_+(x,z) - p_\mu(y,z) \, g_{01}^+(x,y,z)}{p_1(x,y)}
	\nonumber \\
	g^+_{11} &= \frac{1}{p_1(x,y)}\Bigg( \frac{x^d-z^d}{q_+(x,z)} \, g^+_{01}(x,y,z) - \frac{y^d-z^d}{p_\mu(y,z)} \Bigg)
\end{align*}
As above, one verifies that the $g^+_{10}$ and $g^+_{11}$ are indeed polynomials in $x,y,z$.

\medskip

We will now establish that $(g^-,g^+) : \hat P_{a+b+1:\mu-1} \oplus \hat P_{a+b:\mu+1} \longrightarrow \hat P_{a:1} \otimes \hat P_{b:\mu}$ is an isomorphism in $\HMFbi$ (and thereby also in $\HMFbigr$ as $g^\pm$ have $\bC$-degree 0). We do this by employing Remark \ref{rem:H(M)} (see \cite[Cor.\,4.9]{Wu09} for the corresponding graded statement), that is, by showing that $(H(g^-),H(g^+)) : H(\hat P_{a+b+1:\mu-1}) \oplus H(\hat P_{a+b:\mu+1}) \longrightarrow H(\hat P_{a:1} \otimes \hat P_{b:\mu})$ is an isomorphism.

For $\hat P_S$ we have $H(\hat P_\emptyset) = H(\hat P_{\bZ_d}) = 0$ and
$H(\hat P_S) = \bC \oplus \bC$ if $S \neq \emptyset, \bZ_d$. The first case occurs only for $\mu=d-2$, where $H(\hat P_{a+b:\mu+1})=0$.

For $H(\hat P_{a:1} \otimes \hat P_{b:\mu})$ we need to compute the homology of the complex
$$
\xymatrix{
{\begin{array}{c}\bC[y] \\ {\op} \\ \bC[y] \end{array} }
 \ar@/^10pt/[rrrr]^{\left(\begin{array}{rr} \scriptstyle{\eta^{2a+1}y^2} & \scriptstyle{y^{\mu+1}} \\ \scriptstyle{-y^{d-\mu-1}} &  \scriptstyle{-\eta^{-2a-1}y^{d-2}} \end{array}\right)} &&&& 
{\begin{array}{c}\bC[y] \\ {\op} \\ \bC[y] \end{array} }
\ar@/^10pt/[llll]^{\left(\begin{array}{lr} \scriptstyle{-\eta^{-2a-1}y^{d-2}} & \scriptstyle{-y^{\mu+1}} \\ \scriptstyle{y^{d-\mu-1}} &  \scriptstyle{\eta^{2a+1}y^{2}} \end{array}\right)}
}
$$
Define the vectors $v_0 := (\eta^{2a+1} , - y^{d-\mu-3})$ (this has second entry equal to $y^{-1}$ for $\mu=d-2$, but the results below are polynomial nonetheless) and $v_1 = (-y^{\mu-1} , \eta^{2a+1})$. 
One finds
\begin{align*}
	\mathrm{ker}(\bar d_0) &= \bC[y] v_0 \cdot \begin{cases} 1 &; \mu < d-2 \\ y &; \mu = d-2 \end{cases}
	&
	\mathrm{ker}(\bar d_1) &=  \bC[y] v_1
\nonumber	\\
	\mathrm{im}(\bar d_1) &= y^2 \bC[y] v_0
&
	\mathrm{im}(\bar d_0) &= y \bC[y] v_1\cdot \begin{cases} 1 &; \mu < d-2 \\ y &; \mu = d-2 \end{cases}
\end{align*}
Writing $[\cdots]$ for the homology classes, the homology groups 
$H_i := H_i(\hat P_{a:1} \otimes \hat P_{b:\mu})$
are given by
$$
	H_0 = 
	\begin{cases} \{ [v_0] , [yv_0] \} &; \mu < d-2 \\ \{ [yv_0] \} &; \mu = d-2 \end{cases}
	\quad , \qquad
	H_1 = 
	\begin{cases} \{ [v_1] , [yv_1] \} &; \mu < d-2 \\ \{ [v_1] \} &; \mu = d-2 \end{cases}
	\qquad .
$$
The map $(g^-,g^+)$ acts on homology by, for $\mu < d-2$,
\begin{align*}
	H_0(g^-,g^+) &=~
	\Big(~~ \eta^{-2a-1} \beta^- \, [y v_0] ~~,~~ \eta^{-2a-1} [v_0]~~ \Big) \ ,
	\nonumber \\
	H_1(g^-,g^+) &=~ 	\Big( ~~\eta^{-2a-1} [v_1] ~~,~~ \eta^{-2a-1} \beta^+ \, [y v_1]~~ \Big) \ .
\end{align*}
Here $\beta^-$ is the coefficient of $y$ in $g^-_{00}$ and $\beta^+$ is the coefficient of $y$ in $g^+_{01}$.
For $\mu = d-2$, the second entry in the above maps is absent, as $H(\hat P_{a+b:\mu+1})=0$ in this case. Altogether we see that $H(g^-,g^+)$ is indeed an isomorphism.

\medskip

This proves the decomposition 
$\hat P_{a:1} \otimes \hat P_{b:\mu} ~\simeq~
\hat P_{a+b+1:\mu-1} \oplus \hat P_{a+b:\mu+1}$ 
in $\HMFbigr$ and
completes the proof of Theorem \ref{thm:Pdgr-tensor-closed}.

\end{appendices}

\newcommand\arxiv[2]      {\href{http://arXiv.org/abs/#1}{#2}}
\newcommand\doi[2]        {\href{http://dx.doi.org/#1}{#2}}
\newcommand\httpurl[2]    {\href{http://#1}{#2}}

\end{document}